\documentclass[11pt,reqno,final]{article}

\usepackage[color,notref,notcite]{showkeys}
\definecolor{labelkey}{rgb}{0.6,0,1}
\usepackage{epsfig}
\usepackage{a4wide,amssymb,amsmath}
\usepackage{amsthm}
\usepackage{color}
\usepackage{url}
\usepackage{mathtools}

\usepackage[latin1]{inputenc}
\usepackage[english]{babel}

\usepackage{amsfonts}
\usepackage{epsfig}
\usepackage{graphicx}
\usepackage{url}

\usepackage{algorithm}
\usepackage[noend]{algpseudocode}
\floatname{algorithm}{Algorithm}      

\makeatletter
\def\BState{\State\hskip-\ALG@thistlm}
\makeatother

\newcommand{\ba}{\begin{eqnarray}}
\newcommand{\ea}{\end{eqnarray}}

\newtheorem{theo}{\bf Theorem}[section]

\newtheorem{property}[theo]{\bf Property}

\newtheorem{example}[theo]{Example}

\newtheorem{remark}[theo]{\bf Remark}
\providecommand{\keywords}[1]{\textbf{\textit{Keywords: }} #1}

\title{\bf
A High-Order Scheme for Image Segmentation via a modified Level-Set method
}

\author{Maurizio Falcone\footnotemark[1] \and Giulio Paolucci\footnotemark[1] \and Silvia Tozza\footnotemark[2]}
\begin{document}

\maketitle

\renewcommand{\thefootnote}{\fnsymbol{footnote}}

\footnotetext[1]
{
 Dipartimento di Matematica,  ``Sapienza" Universit{\`a}  di Roma,
 P.le Aldo Moro, 5 - 00185 Rome, Italy
 ({\tt e-mail:falcone@mat.uniroma1.it, paolucci@mat.uniroma1.it})
}


\footnotetext[2]
{
Istituto Nazionale di Alta Matematica, U.O. Dipartimento di Matematica,  ``Sapienza" Universit{\`a}  di Roma,
P.le Aldo Moro, 5 - 00185 Rome, Italy
({\tt e-mail: tozza@mat.uniroma1.it})  \\
The three authors are members of the INdAM Research group GNCS. \\
}

\renewcommand{\thefootnote}{\arabic{footnote}}

\abstract{In this paper we propose a high-order accurate  scheme for image segmentation based on the level-set method. In this approach, the curve evolution is described as the 0-level set of a representation function but we modify the velocity that drives the curve to the boundary of the object in order to obtain a new velocity with additional properties that are extremely useful to develop a more stable high-order approximation with a small additional  cost. The approximation scheme proposed here is the first  2D version of an adaptive ``filtered" scheme recently introduced and analyzed  by the authors in 1D. This approach is interesting since the implementation of the filtered scheme is rather efficient  and easy. The scheme combines two building blocks (a monotone scheme and a high-order scheme) via a filter function and smoothness indicators that allow to detect the regularity of the approximate solution adapting the scheme in an automatic way. Some numerical tests on synthetic and real images confirm the accuracy of the proposed method and the advantages given by the new velocity.  
}
\vskip12pt


\section{Introduction}\label{sec:introduction}
The Level-Set (LS) method  has been introduced by Osher and Sethian in the 1980s \cite{Sethian85,OS88} and then used to deal with several applications, e.g.  fronts propagation, computer vision, computational fluids dynamics (see the monographs by Sethian \cite{Sethian99} and by Osher and Fedkiw \cite{OF03} for several interesting examples). This method is nowadays very popular for its simplicity and for its capability to deal with topological changes. In fact, the main advantage of the LS method is the possibility to easily describe time-varying objects, follow shapes that change topology, for example when a shape splits in two, develops holes, or the reverse of these operations. For the image segmentation problem the application of the LS method is based on the evolution of a curve according to a 
normal velocity based on the gray-levels of the image, typically the curve is described as the 0-level set of a representation function (or level-set function).\\
In $\mathbb R^2$ the LS method corresponds to define  an initial  closed curve $\Gamma_0$  using an auxiliary function $v_0$ which has to change sign 
on $\Gamma_0$. The evolution of that curve at time $t$ is denoted by $\Gamma_t$ and  is represented by the 0-level set of a function $v$, i.e.
\begin{equation}\label{def:gammat}
\Gamma_t:=\{(x,y):v(t,x,y)=0\}.
\end{equation}
This function $v$ is the unique viscosity solution of the following evolutive nonlinear equation of Hamilton-Jacobi type
\begin{equation}\label{eq:HJ_eikonal} 
\left\{
\begin{array}{ll}
v_t+c(t,x,y)|\nabla v|=0,\qquad&(t,x,y)\in \mathbb (0,T)\times \mathbb R^2,\\
v(0,x,y)=v_0(x,y), &(x,y)\in \mathbb R^2,
\end{array}
\right.
\end{equation}
where $\nabla v := (v_x,v_y)$ denotes the spatial gradient of $v$. Usually, the velocity $c(t, x,y)$  does not change sign during the evolution and the orientation depends on the type of evolution (outward for an expansion and inward for a shrinking). 
Typically $v_0$ must be a \emph{proper representation} of the initial front $\Gamma_0$, satisfying
\begin{equation}
\left\{
\begin{array}{ll}
v_0(x,y)<0,\qquad&(x,y)\in \Omega_0,\\
v_0(x,y)=0,\qquad&(x,y)\in\Gamma_0,\\
v_0(x,y)>0,\qquad&(x,y)\in \mathbb R^2\setminus\overline\Omega_0,
\end{array}
\right.
\end{equation}
where $\Omega_0$ is the region delimited by $\Gamma_0$ (or the reverse inequalities). \\
The LS method can handle velocities also depending on physical quantities in order to describe several phenomena. Typical examples are:
\begin{enumerate}
	\item[a)] $c(t, x,y)$, isotropic growth with time varying velocity
	\item[b)]  $c(t,x,y,\eta)$, anisotropic growth, dependent on normal direction
	\item [c)] $c(t,x,y,k)$, Mean Curvature Motion, with $k(t,x,y)$ mean curvature to the front at time $t$. 
\end{enumerate}
The literature on the LS method is huge as well as the range of applications where the method has been successfully applied. We refer the interested reader to \cite{CCCD93,MSV95} for the isotropic case, to \cite{KKOTY96} for the anisotropic growth and to \cite{CKS97, CFF10} for the curvature case (see also the monographs \cite{Sethian99, OF03} and the references therein).  
Depending on the choice of the velocity,  the front evolution will be described by first or second order partial differential equations, e.g. case c) gives rise to a degenerate second order equation.  Here we limit ourself to first order problems \eqref{eq:HJ_eikonal} and velocities which depend only on time and position and we consider a velocity sign constant in time since this is the more relevant case for image segmentation (a change in sign is relevant for dislocations in material science and has been studied in \cite{ CCF06,CFFM08}). \\
To set our contribution into perspective, let us mention that the segmentation problem has been solved by various techniques which mainly rely on two different approaches: variational methods and active contour methods. For the first approach the interested reader can look at \cite{CCZ13,BW04} and the references therein. 
For the link between the two classes of methods, see \cite{CKS97}.  To have an idea of other segmentation methods, we refer the reader to the surveys \cite{ZMCL16, CFF19} and to the papers \cite{CC92,ZZSZ10} for the so called ``balloon model"  first introduced in \cite{CC91}. This method is based on the introduction of a potential giving a driving force to the segmentation process and typically leads to second order partial differential equations that are at present outside the reach of the filtered scheme proposed in this paper. For high-order Runge-Kutta methods in the framework of image segmentation we refer to \cite{SR09}. 
As already said, we will apply the level set method based on \eqref{eq:HJ_eikonal} looking for an accurate numerical method. High-order methods have been proposed for \eqref{eq:HJ_eikonal} and most of them are based on non oscillatory local interpolation techniques that allow to avoid spurious oscillations around discontinuities of the solution and/or jumps in the gradient. 
These techniques were originally developed for conservation laws (see the seminal paper \cite{HEOC86} and the references therein), the research activity on essentially non oscillatory (ENO)  methods has been rather effective and a number of improvements have been proposed in e.g. \cite{LOC94, JS96, AB03} so that now ENO and Weighted ENO (WENO) techniques are rather popular in many applications (see \cite{ZS16} for a recent survey). We should also mention that later on   these techniques were successfully applied  to the numerical solution of  Hamilton-Jacobi equations \cite{OS91} opening the way to other applications. The reader  interested in image processing will find an application to image compression in  \cite{AACD02,ABM10} and  an application to image segmentation in \cite{ZZLZ16}. However, a general convergence theorem for ENO/WENO schemes is still missing and their application is a rather delicate issue.  These limitations have motivated further investigations and  a new class of high-order methods for evolutive Hamilton-Jacobi equations have been proposed looking for a different approach based on {\em ``filtered schemes" } introduced in this framework by Lions and Souganidis  \cite{LS95} (it is important to note that the name does not refer to filtering a noise as is common in the imaging community but rather to the presence of a filter function as we will explain later). The class of filtered schemes is based on a simple  coupling of  a monotone scheme with a high-order scheme. Monotone schemes are convergent to the weak (viscosity) solution but they are  known to be at most first order accurate, whereas  high-order schemes gives a higher accuracy but  in general are not stable. The crucial point is the coupling between the two schemes which is obtained via a {\em filter function} that selects which scheme has to be applied at a  node of the grid  in order to guarantee  (under appropriate assumptions) a global convergence. The construction of these schemes is rather simple as explained by Oberman and Salvador \cite{OS15} because one can couple various numerical methods and leave the filter function deciding the  switch between the two schemes. A general convergence result has been proved by Bokanowski, Falcone and Sahu in \cite{BFS16} and recently improved by Falcone, Paolucci and Tozza \cite{FPT18a} with an adaptive and automatic choice of the parameter governing the switch. Note that the adaptation of the parameter depends on some regularity indicators in every cell, these indicators are computed at each iteration and this guarantees convergence. Some contributions to extend  filtered schemes to second order problems can be found in \cite{FO13} for the stationary case and in \cite{BPR18} for the evolutive case and financial applications. 
\paragraph{Our contribution}
The contribution of this paper is twofold: from a theoretical point of view, we propose a modified velocity for the  segmentation problem in the level-set approach. This new velocity is important (and necessary) for the numerical approximation of high-order schemes, like the filtered scheme considered here and its properties will allow to avoid the re-initialization procedure used e.g. in \cite{MSV95}. We also improve the accuracy of the method by applying an 
adaptive high-order filtered scheme to the segmentation problem. 
This requires a 2D extension of the scheme proposed in 1D in \cite{FPT19}, for which convergence has been proved under rather general assumptions in \cite{FPT18a}, thanks to the definition of new full 2D smoothness coefficients (see Sect. \ref{sec:AFS}).  
We will show that this choice  is competitive with respect to other high-order schemes, as the WENO scheme, in terms of accuracy and computational cost, analyzing several experiments. This contribution is more efficient in terms of CPU time  and is interesting from a theoretical point of view since a precise convergence result for WENO scheme is still missing. 
\paragraph{Paper organization}
The paper is organized as follows: In Sect. \ref{sec:image_segm_via_LS}, we give the idea behind the new definition of the velocity function, with details on its construction and explaining why it is so important to introduce it when applying high-order approximation schemes. In Sect. \ref{sec:AFS}, we briefly present the Adaptive Filtered (AF) scheme and recall some of its basic elements. In Sect. \ref{sec:num_impl}, we give some information on the implementation and we give a sketch of the \emph{Algorithm} \ref{pseudo_code} for the solution of the segmentation problem via the new AF scheme. Finally, the performances of this new method for the segmentation problem are illustrated in Sect. \ref{sec:tests} where we compare it with a classical  monotone  scheme (first order accurate) and another high-order scheme (the WENO scheme) on a series of virtual and real images presenting a detailed error analysis of the numerical experiments. 

\section{Image segmentation via a modified LS method}\label{sec:image_segm_via_LS}
The boundaries of one (or more) object(s) inside a given image are characterized by an abrupt change of the intensity values $I(x,y)$ of the image, so that the magnitude of $|\nabla I(x,y)|$ can be used as an indication of the edges. Let us assume that $\nabla I $ exists at least almost everywhere, the definition of the velocity $c(x,y)$ in \eqref{eq:HJ_eikonal} plays a crucial role and must be defined in a proper way in order to guarantee that the curve evolution is close to $0$ when the front is close to an edge, this will stop the evolution.  The normal velocity will be positive or negative depending on the case, expanding or shrinking, respectively, and if the velocity is just given in term of $I$ it will not change sign during the evolution, as it could happen when the velocity depends on other geometrical properties of the curve (e.g. the curvature). We will focus on velocities which ignore the curvature since the numerical approach we present here is well adapted to first order evolutive Hamilton-Jacobi equations (the extension to second order problems goes beyond the scopes of this paper). 
Note that, in order to reduce the noise, several methods have been proposed in the literature. A very simple one is to use the convolution with a Gaussian kernel. This can be obtained by evolving the original function $I$ of the gray levels according to the heat equation for a short time interval (see Sect. \ref{sec:num_impl} for more details). Due to the regularizing effect of the heat equation this also guarantees that $\nabla I$ exists.\\
\noindent Several definitions of the velocity function $c(x,y)$ have been proposed in literature. A typical example is
\begin{equation}\label{vel_c1}
c_1(x,y)=\frac{1}{\left(1+|\nabla(G\ast I(x,y))|^\mu\right)},\qquad \mu\geq 1,
\end{equation}
where $\mu$ is used to give more weight to the changes in the gradient, if necessary. 
In \cite{CCCD93} the authors proposed that velocity with $\mu=2$, and in  \cite{MSV95} with $\mu=1$. 
According to this definition, the velocity takes values in $[0,1]$ and   has values that are close to zero at points where  the image gradient is high and equal to 1  where  $I$ is constant. \\
Another possible choice has been proposed in  \cite{MSV95} and has the form
\begin{equation}\label{vel_c2}
c_2(x,y)=1-\frac{|\nabla(G\ast I(x,y))|-m}{M-m},
\end{equation}
where  $m$ and $M$ are respectively the minimum and the maximum values of $|\nabla(G\ast I(x,y))|$. This  velocity has similar properties with respect to \eqref{vel_c1} and takes values in $[0,1]$ but is close to $0$ if the magnitude of the image gradient is close to its maximal value, and equal to $1$ otherwise. \\
It is clear that both definitions have the desired properties, but with slightly different features. More precisely, in the first case the velocity depends more heavily on the changes in the magnitude of the gradient, allowing for an easier detection of the edges but  possibly producing false edges inside the object (e.g. when specular effects are present in the image).
The velocity \eqref{vel_c2} is smoother inside the objects, being less dependent on the relative changes in the gradient, but might present some problems in the detection of all the edges if at least one of those is ``more marked''. 

\subsection{Extension of the velocity function}\label{subsec:extension_of_c}
The edge-stopping function which is defined choosing one of the above mentioned velocities,  has a physical meaning only on the front $\Gamma_t$ since it was designed precisely to force the $0$-level set to stop close to the edges. As it has been observed  in \cite{MSV95},  its meaning does not come from the geometry of $v$ but only from the configuration of the front $\Gamma_t$. 
Using one of the  classical velocity functions introduced before, as it will be clarified by the numerical tests in Sect. \ref{sec:tests}, high-order schemes produce unstable results since the numerical approximation of $v$ can start to produce spurious oscillations near the edges where the front should stop. 
This problem can be solved by adding a limiter as in \cite{BFS16} but here we present a different technique that avoids the use of a limiter and adapts automatically the scheme according to the regularity of the solution in order to produce more accurate  results. To this end, we need to extend the image-based velocity function $c(x,y)$ to all the level sets of the representation $v$ in order to give a physical meaning  to the speed used in the whole domain. 

Following some of the ideas discussed in \cite{MSV95}, we want to extend in a simple way the velocity to the whole domain. Our approach exploits the choice of the initial condition $v_0$ (which is free)  and allows to avoid all the heavy computations required by the numerical approach solution proposed by the authors in \cite{MSV95}. The modification has an interesting interpretation in terms of the method of characteristics, as we will see later in this section, and  allows for stable numerical results. Thus, recalling their approach, the first property that the velocity has to satisfy is: 
\begin{property}
	Any external (image-based) speed function that is used in the equation of motion written for the function $v$ should not cause the level sets to collide and cross each other during the evolutionary process.
\end{property}
To present the main idea, let us consider the signed distance to the initial $0$-level set as the initial representation function. This is the  classical choice
\begin{equation}\label{eq:v0_distance}
v_0(x,y)=
\left\{
\begin{array}{ll}
- dist\left\{(x,y),\Gamma_0\right\} \hbox{ in } \Omega_0\\
dist\left\{(x,y),\Gamma_0\right\} \hbox { outside } \Omega_0
\end{array}
\right.
\end{equation}
where $\Omega_0$ is the internal region delimited by $\Gamma_0$. 
Therefore, with this choice we can define the velocity extension as follows: 
\begin{property}\label{propertyLS}
	The value of the speed function $c(x,y)$ at a point $P$ lying on a level set $\{v=C\}$ is exactly the value of $c(x,y)$ at a point $Q$, such that the point $Q$ is a distance $C$ away from $P$ and lies on the level set $\{v=0\}$.
\end{property}
\noindent Note that the point $Q$ is uniquely determined whenever the normal direction in $P$ is well defined. In fact, $Q = P - c(x,y) \eta(P)$, where $\eta$ is the outgoing normal and this will provide a definition also when the level sets are non convex.  To develop a formal argument we will assume that the normal is sufficiently regular.\\
In order to compute the point $Q$ on the $0$-level set associated to each point $P$ of any level set, the authors  introduce in \cite{MSV95} on pages 162-164 a neighborhood of the 0-level set called \emph{Narrow-Band}. They fix the width $\delta$ of the band  and the related number of iterations $l$ needed for the $0$-level set to reach the boundary of the narrow-band (this is done by the procedure explained at page 164 of [20]). This procedure forces to modify the representation function since the solution is updated only inside the band and a reinitialization step becomes strictly necessary in order to restore the meaning of the distance function. Quoting from their paper, such procedure either requires at least $O(N^3)$ computations for a grid of $N$ points in each direction or requires to compute  the solution of an associated stationary eikonal equation. 
In order to reduce the computational complexity, the method proposed here is based on a direct assignment of the associated point on the $0$-level set, which requires only $O(N^2)$ operations to compute the modified velocity at each iteration, and greatly simplifies the problem. Moreover, it allows for the use of representation functions (i.e. initial conditions) even more regular than the signed distance function. 
The idea is straightforward and is based on the fact that the evolution is oriented in the normal direction to the front. If the reciprocal position of the level sets is also known (that is why we must choose wisely the initial condition) and we make all the points in the normal direction to the $0$-level set evolve according to the same law  it is reasonable to expect that all such points will keep their relative distance unchanged as time flows. \\
To illustrate and motivate our modification, let us still consider the distance to $\Gamma_0$ \eqref{eq:v0_distance} as initial condition and let us consider the shrinking case as example. Then, by construction, all the $C$-level sets are at a distance $C$ from the $0$-level set, as stated by Property \ref{propertyLS}. Hence, if we consider a generic point $(x_c,y_c)$ on a $C$-level set, then it is reasonable to assume that the closest point on $\Gamma_0$ should be
\begin{equation}
(x_0,y_0)=(x_c,y_c)-v(t,x_c,y_c)\frac{\nabla v(t,x_c,y_c)}{|\nabla v(t,x_c,y_c)|}.
\end{equation}
Therefore, it seems natural to define the extended velocity $\widetilde c(x,y)$ as
\begin{equation}
\widetilde c(x,y,v,v_x,v_y)=c\left(x-v\frac{v_x}{|\nabla v|},y-v\frac{v_y}{|\nabla v|}\right),
\end{equation}
which coincides with $c(x,y)$ on the $0$-level set, as it is needed. The same approach can be applied as long as the initial distance between the level sets is known. In that case, if we want higher regularity to the evolving surface, which would be preferable in the case of high-order schemes such as  those we use in the numerical tests, we can define an appropriate initial condition, for example, by simply rotating a regular function in one space dimension. More precisely, let us consider a regular function $\overline{v}_0:\mathbb R^+\to \mathbb R$ such that $\overline{v}_0(r_0)=0$, where $r_0$ is the radius of the initial circle $\Gamma_0$ (e.g. the right branch of a parabola centered in the origin), and let us define $v_0(x,y)$ rotating its profile, that is
\begin{equation}
\label{u0_simm_circ}
v_0(x,y)=\overline v_0\left(\sqrt{x^2+y^2}\right).
\end{equation}
Then, it is clear that the C-level sets of $v_0$ are located at a distance
\begin{equation}
\label{distance_gen}
d(C):=\overline{v}_0^{-1}(C)-r_0,\qquad\textrm{with }\overline v_0^{-1}(C)\geq 0,
\end{equation}
from the $0$-level set and, according to our previous remarks, they should keep this property as time evolves. Consequently, also in this case we can define
\begin{equation}\label{eq_new_vel}
\widetilde c(x,y,v,v_x,v_y)=c\left(x-d(v)\frac{v_x}{|\nabla v|},y-d(v)\frac{v_y}{|\nabla v|}\right).
\end{equation}
More details on the function $d(v)$ will be given in Sect. \ref{subsec:motivation_ctilde}. 
For simplicity, in the last construction we assumed the representation function to be centered in the origin, but it is straightforward to extend the same procedure to more general situations. Note also that if we have only one object to be segmented (or we are considering the shrinking from the frame of the picture) we can always use a representation function centered in the origin since we can choose freely  the domain of integration, given by the pixels of the image. 

\subsection{Motivations of the new velocity function}\label{subsec:motivation_ctilde}
The modification of the velocity $c(x,y)$ into $\widetilde c(x,y,v,v_x,v_y)$ defined in  \eqref{eq_new_vel} with $d(v)=0$ if $v=0$, is to follow the evolution of the $0$-level set and then to define the evolution on the other level sets accordingly.  This allows to give a geometrical interpretation of the new velocity and to establish some properties that will also guarantee existence and uniqueness for the first order evolutive problem. 
The new velocity can be seen as a \emph{characteristic based velocity}. As a first step let us analyze the characteristics of the equation, in these computations we  assume that we always have the necessary regularity.  In particular, we assume $v\in C^2(\Omega)$ (or at least $C^2$ in space and $C^1$ in time) and $c(x,y)\in C^1(\Omega)$. We introduce the notations of the vectors $z := (x,y)$ and $p:= (p_1,p_2) = (v_x,v_y)$ that will be used only in this section. Using these notations, the Hamiltonian in our case can be written as
\begin{equation}
H(z,v,p)=\widetilde c(z,v,p)|p|. 
\end{equation}
Let us introduce the method of characteristics, writing the usual system
\begin{equation}
\label{eq:system_char_segm}
\left\{
\begin{array}{l}
\dot{z}(s)=\nabla_pH\\
\dot{v}(s)=\nabla_p H \cdot p-H\\
\dot{p}(s)=-\nabla H - H_v  p, 
\end{array}
\right.
\end{equation}
where $\dot{f}$, $f=z,v,p$  denotes the derivative with respect to the variable $s$, $\nabla_p H$ the gradient with respect to $p$, $H_v$ the partial derivative with respect to the (scalar) value $v$, and $\nabla H$ the usual spatial gradient of $H$ already introduced in Sect. \ref{sec:introduction}. 
In our case, defining for brevity the point $ (\xi, \zeta) := (x-d(v)\frac{p_1}{|p|},y-d(v)\frac{p_2}{|p|})$,  
we obtain
\begin{align}
\frac{\partial H}{\partial p_1}=& \frac{\partial \widetilde c}{\partial p_1} |p| + \widetilde c \frac{\partial |p|}{\partial p_1} \\
=& \left( \frac{\partial c}{\partial \xi} \cdot \frac{\partial \xi}{\partial p_1} + \frac{\partial c}{\partial \zeta} \cdot \frac{\partial \zeta}{\partial p_1} \right) |p| + \widetilde c  \frac{ p_1}{ |p|}\nonumber\\
=&-d(v) \frac{\partial c}{\partial \xi} \left(\frac{|p|-\frac{p^2_1}{|p|}}{{|p|}^2}\right)|p| -d(v) \frac{\partial c}{\partial \zeta} \left(-\frac{p_1 p_2}{|p|^3}\right)|p|+\widetilde c \frac{p_1}{|p|} \nonumber\\
=&-d(v)  \frac{\partial c}{\partial \xi} \frac{p_2^2}{{|p|}^2}+d(v)  \frac{\partial c}{\partial \zeta} \frac{p_1 p_2}{{|p|}^2}+\widetilde c \frac{p_1}{|p|}, \nonumber
\end{align}
and analogously $\frac{\partial H}{\partial p_2}$, so that we obtain
\begin{equation}
\label{def_Hp}
\nabla_p H=
\left(
\begin{array}{l}
\frac{d(v) p_2}{|p|^2}\left(p_1 \frac{\partial c}{\partial \zeta} - p_2  \frac{\partial c}{\partial \xi} \right)+\widetilde c \frac{p_1}{|p|} \\
\frac{d(v) p_1}{|p|^2}\left(p_2 \frac{\partial c}{\partial \xi} - p_1  \frac{\partial c}{\partial \zeta} \right)+\widetilde c \frac{p_2}{|p|} \\
\end{array}
\right)
\hbox{ that implies }
\nabla_p H \cdot p=\widetilde c(z,v,p)|p|. 
\end{equation}
Therefore, the system \eqref{eq:system_char_segm} becomes in our case the following: 
\begin{equation}
\left\{
\begin{array}{l}
\dot{z}(s)=\nabla_p H\\
\dot{v}(s)=\widetilde c(z,v,p)|p|-\widetilde c(z,v,p)|p|=0\\
\dot{p}(s)=-\nabla \widetilde c(z,v,p) |p|+d^\prime(v)\nabla \widetilde c(z,v,p)|p|^2=\nabla \widetilde c(z,v,p)|p|(d^\prime(v)|p|-1), 
\end{array}
\right.
\end{equation}
where $d^\prime(v)$ denotes the derivative of the distance $d$ with respect to the value $v$. \\
Now, choosing $d^\prime(v)$ such that
\begin{equation}
\label{eq_dpu}
d^\prime(v)=|p|^{-1},
\end{equation}
we have the final system
\begin{equation}
\label{final_system}
\left\{
\begin{array}{l}
\dot{z}(s)=\nabla_p H, \\
\dot{v}(s)=0, \\
\dot{p}(s)=0, 
\end{array}
\right.
\end{equation}
which states that, as long as the function $\widetilde c$ remains smooth enough ($ \frac{\partial c}{\partial \xi} \approx 0$ and $ \frac{\partial c}{\partial \zeta} \approx 0$), the characteristics are basically directed in the normal direction and  along them both the height and the gradient are preserved. 
Looking at the third relation of \eqref{final_system} and at the choice  \eqref{eq_dpu}, since $p(s)\equiv p(0)=\nabla v_0$ 
along the characteristics, we can choose  simply
\begin{equation}\label{eq_dpu_2}
d^\prime(v)=|\nabla v_0|^{-1},
\end{equation}
which is the trivial case with the function $d(v)=v$ and also for $d(v)$ given by the previous definition \eqref{distance_gen}. 
In fact, using the inverse function theorem, we have
\begin{equation}
d^\prime(v)=\frac{d}{d v}\left(\overline v_0^{-1}(v)\right)=\frac{1}{\overline v^\prime_0(w)},
\end{equation}
with $w$ such that $\overline v_0(w)=v$. Moreover, recalling the definition (\ref{u0_simm_circ}), we can compute
\begin{align}
|\nabla v_0(x,y)|&=\left|\nabla \overline v_0\left(\sqrt{x^2+y^2}\right)\right|\\
&=\left|\left(\frac{\overline v_0^\prime \left(\sqrt{x^2+y^2}\right)x}{\sqrt{x^2+y^2}},\frac{\overline v_0^\prime \left(\sqrt{x^2+y^2}\right)y}{\sqrt{x^2+y^2}}\right)\right| \nonumber \\
&=\frac{\overline v_0^\prime \left(\sqrt{x^2+y^2}\right)}{\sqrt{x^2+y^2}} \sqrt{x^2+y^2}=\overline v_0^\prime \left(\sqrt{x^2+y^2}\right), \nonumber
\end{align}
and then it is enough to consider $(x,y)$ such that $w=\sqrt{x^2+y^2}$. 
From a numerical point of view, we compute the function $d(v)$ analytically, by using \eqref{eq_dpu_2} that exploits our knowledge of the initial condition $v_0$, e.g. given by \eqref{u0_simm_circ}. 
Thanks to the previous computations, we reached a good understanding of the nature of the evolution given by (\ref{eq:HJ_eikonal})-(\ref{eq_new_vel}), but we still have not justified the main motivation that led us to define (\ref{eq_new_vel}), that is to make all the level sets of $v$ evolve according to the same law. More precisely, we have to show that, if we consider the evolution of two points on the same characteristic but on two different level sets, say the $0$-level set $z^0(s)$ and a generic level set $z^\ell(s)$, then their relative distance (along the characteristic) does not change during the evolution. This fact would imply that, if we choose the level sets of $v_0$ to be such that 
\begin{equation}\label{init_cond_x}
z^0(0)=z^\ell(0)-d(v_0(z^\ell))\frac{\nabla v_0(z^\ell)}{|\nabla v_0(z^\ell)|},
\end{equation}
then the points $\underline z(s):=z^\ell(s)-d(v(z^\ell))\frac{p(z^\ell)}{|p(z^l)|}$ are always on the $0$-level set of $v$. 
In order to prove this last statement, let us proceed by a simple differentiation, dropping the dependence on $z^\ell$ for brevity,
\begin{align}
\dot {\underline z}(s)=& \dot z^\ell(s)-\frac{d}{ds}\left(d(v)\frac{p}{|p|}\right)  \\
=&\dot z^\ell(s) -\left[d^\prime(v)\dot v(s)\frac{p}{|p|}+\frac{d(v)}{|p|^2}\left(\dot p(s)|p|-\frac{d}{ds}(|p(s)|)p\right)\right]. \nonumber
\end{align}
Recalling the relations in the system \eqref{final_system} with respect to $\dot{v}(s)$ and $\dot{p}(s)$, we can write
\begin{align}
\dot {\underline z}(s)=&\dot z^\ell(s)+\frac{d(v)}{|p|^2}\left(\frac{p\cdot \dot p(s)}{|p|}\right)p \\ 
=&\dot z^\ell(s). \nonumber
\end{align}
This calculation shows  that the points  $\underline z(s)$ and $z^\ell(s)$ evolve according to the same law along characteristics. 
Note that if \eqref{init_cond_x} holds then $\underline z(s)\equiv z^0(s)$ till the characteristics do not cross. In fact, computing the total derivative with respect to $s$, we have
\begin{equation}
\frac{d}{ds} v(s,\underline z(s))=v_s + v_x \dot{\underline z}(s)=v_s + v_x \dot{z}^\ell(s)=\frac{d}{ds} v(s, z^\ell(s))=0,
\end{equation}
since the points of $z^\ell(s)$ are on the same level set, and so are those of $\underline z(s)$, as we wanted. 
This directly implies that the points $\underline z:= \left(z-d(v)\frac{\nabla v}{|\nabla v|}\right)$ are on the $0$-level set of $v$ as long as the gradient is preserved. \\
In order to reduce the computational cost and to simplify the implementation, as will be better explained in Sect.~\ref{sec:num_impl}, we decided to split the computation of the Hamiltonian into two steps. First we compute the points $\underline z$ and then we update the numerical solution of the full problem (\ref{eq:HJ_eikonal}) with the velocity (\ref{eq_new_vel}) via the following \emph{simplified problem} with isotropic velocity 
\begin{equation}\label{eq_simplified}
v_t+c(\underline z)|\nabla v|=0,\qquad(t,x,y)\in \mathbb (t_n,t_{n+1})\times \mathbb R^2,
\end{equation}
where, as usual, $t_n=t_0+n\Delta t$, and $\Delta t$ is the time step. 
Using this procedure the dependencies on the gradient and on the value of $v$ are frozen at every iteration, leaving just an explicit dependence on the variable $t$. This is why we consider velocities depending only on space variables when defining the numerical schemes. \\
Let the original velocity $c(x,y)$  be Lipschitz continuous (this is the classical assumption) and let us denote by $L_c$ its Lipschitz constant. We conclude showing that  our modified velocity $\tilde c$ is still Lipschitz continuous. This point is important to guarantee existence and uniqueness of the characteristics and of the viscosity solution for the evolutive problem \eqref{eq:HJ_eikonal} driven by the new velocity. Let $z=(x,y)$ and $z'=(x',y')$ be two points in the plane and $\underline z$ and $\underline z'$ the corresponding points used in the definition of the new velocity $\tilde c$. Provided that the solution and the normal vector are Lipschitz continuous, we have
\begin{eqnarray}
|c(\underline z)-c(\underline z')|&&= \left|c\left (z-d(v(z))\frac{\nabla v(z)}{|\nabla v(z)|}\right)-c\left(z'-d(v(z'))\frac{\nabla v(z')}{|\nabla v(z')|}\right)\right|\nonumber \\
&& \le L_c \left( |z-z'|+ |d(v(z')) \eta (z') -d(v(z)) \eta (z)|\right)\\
&&\le L_c \left( |z-z'|+ |d(v(z'))  -d(v(z))| |\eta (z')| + |d(v(z))|| \eta (z')- \eta(z)| \right)\nonumber\\
&&\le L_c(1+C_1 + C_2) |z' -z| \nonumber
\end{eqnarray}
where $\eta$ represents the normal unitary vector at the point and $C_1$, $C_2$ are two appropriate constants (remember that the distance from the level set stays  bounded during the evolution). 

\section{The Adaptive Filtered Scheme}\label{sec:AFS}
In this section we will introduce and illustrate the AF scheme we will use to approximate the viscosity solution of the problem \eqref{eq:HJ_eikonal}. 
It is important to note that the name does not refer to filtering a noise as is common in the imaging community but rather to the presence of a filter function as we will explain later in this section. 
For more details on the AF scheme, see \cite{FPT19}. 
We assume that the Hamiltonian $H$ and the initial data $v_0$ are Lipschitz continuous functions in order to ensure the existence and uniqueness of the viscosity solution \cite{CL83}. 
For a detailed presentation of uniqueness and existence results for viscosity solutions, we refer the reader to \cite{CL83} and \cite{Barles94}.

Now, let us define a uniform grid in space $(x_j,y_i)=(j\Delta x,i\Delta y)$, $j$,$i \in \mathbb Z$, and in time $t_n=t_0+n\Delta t$, $n\in[0,N_T]$, with $(N_T-1)\Delta t< T\leq N_T\Delta t$. Then, we compute the numerical approximation $u^{n+1}_{i,j}=u(t_{n+1},x_j,y_i)$ with the simple formula
\begin{equation}
\label{AFS_2D}
u_{i,j}^{n+1}=S^{AF}(u^n)_{i,j}:=S^M(u^n)_{i,j}+\phi_{i,j}^n\varepsilon^n\Delta t F\left(\frac{S^A(u^n)_{i,j}-S^M(u^n)_{i,j}}{\varepsilon^n \Delta t}\right),
\end{equation}
where  $S^M$ and $S^A$ are respectively the monotone and the high-order scheme dependent on both space variables,  
$F: \mathbb{R} \to \mathbb{R}$ is the \emph{filter function} needed to switch between the two schemes, 
$\varepsilon^n$ is the switching parameter at time $t_n$, and $\phi_{i,j}^n$ is the \emph{smoothness indicator function} at the node $(x_j,y_i)$ and time $t_n$, based on the \emph{2D-smoothness coefficients} defined in \cite{PaolucciPhD} and briefly recalled later on in this section. 
The AF scheme here introduced is convergent, as proven in \cite{FPT18a}. 
In the sequel the gradient components will be denoted by the usual notation $(p,q)$ and $p^+$, $p^-$ will be right and left discrete derivatives with respect to $x$ (similar notations apply to $q$ that denotes the discrete partial derivative with respect to $y$). 

The two schemes composing the AF scheme can be freely chosen, provided that they satisfy the following assumptions:\\
\noindent
\textbf{Assumptions on $S^M$:} 
The scheme is consistent, monotone and can be written in \emph{differenced form}
\begin{equation}\label{eq:FD_2D}
u^{n+1}_{i,j} =  S^{M}(u^{n})_{i,j} := u^{n}_{i,j} -\Delta t ~ h^M\left(x_j,y_i,D_x^-u^n_{i,j},D_x^+u^n_{i,j},D_y^-u^n_{i,j},D_y^+u^n_{i,j}\right)
\end{equation}
for a Lipschitz continuous function $h^M(x,y,p^-,p^+,q^-,q^+)$, with $D_x^{\pm} u^n_{i,j}:=\pm \frac{u^n_{i,j\pm 1}-u^n_{i,j}}{\Delta x}$ and $D_y^{\pm}u^n_{i,j}:=\pm \frac{u^n_{i\pm 1,j}-u^n_{i,j}}{\Delta y}$. \\
\noindent
\textbf{Assumptions on $S^A$:} The scheme has a high-order consistency and can be written in \emph{differenced form}
\begin{align}
\label{eq:HA_2D}
u^{n+1}_{i,j} = S^{A}(u^{n})_{i,j} :=u^{n}_{i,j}-\Delta t  h^{A}&\left(x_j,y_i,D_{k,x}^{-}u_{i,j},\dots, D_x^-u^n_{i,j},D_x^+u^n_{i,j},\dots,D_{k,x}^{+}u^n_{i,j},\right.\nonumber\\
&\left.\quad D_{k,y}^{-}u_{i,j},\dots, D_y^-u^n_{i,j},D_y^+u^n_{i,j},\dots,D_{k,y}^{+}u^n_{i,j}\right),
\end{align}
for a Lipschitz continuous function $h^A(x,y,p^-,p^+,q^-,q^+)$ (in short), with 
\begin{equation*}
D^{\pm}_{k,x} u^n_{i,j}:=\pm \frac{u^n_{i,j\pm k}-u^n_{i,j}}{k\Delta x} \qquad \textrm{and} \qquad D^{\pm}_{k,y} u^n_{i,j}:=\pm \frac{u^n_{i\pm k,j}-u^n_{i,j}}{k\Delta y}.
\end{equation*}

\begin{example} 
	As examples of monotone schemes in differenced form satisfying the hypotheses stated before, we can consider  the simple numerical hamiltonian
	\begin{equation}
	\label{hm_eik_2D}
	h^M(p^-,p^+,q^-,q^+):=\sqrt{{\max\{p^-,-p^+,0\}}^2+{\max\{q^-,-q^+,0\}}^2}
	\end{equation}
	for the \emph{eikonal equation}
	\begin{equation}
	v_t+\sqrt{v_x^2+v_y^2}=0,
	\end{equation}
	or, for more general equations also depending on the space variables, we can use the 2D-version of the \emph{local Lax-Friedrichs hamiltonian}
	\begin{align}
	\label{local_lax_fried_2D}
	h^M(x,y,p^-,p^+,q^-,q^+) := &H\left(x,y,\frac{p^++p^-}{2},\frac{q^++q^-}{2}\right)\nonumber\\
	&-\frac{\alpha_x(p^-,p^+)}{2}(p^+-p^-)-\frac{\alpha_y(q^-,q^+)}{2}(q^+-q^-),
	\end{align}
	with
	\begin{equation}
		\alpha_x(p^-,p^+) := \max_{\genfrac{}{}{0pt}{2}{x,y,q,}{p \in I(p^-,p^+)}}\left|H_p(x,y,p,q)\right|,\qquad \alpha_y(q^-,q^+) := \max_{\genfrac{}{}{0pt}{2}{x,y,p,}{q \in I(q^-,q^+)}}\left|H_q(x,y,p,q)\right|,
	\end{equation}
	where $I(a,b) := [\min(a,b),\max(a,b)]$. 
	This scheme is monotone under the restrictions $\frac{\Delta t}{\Delta x}\cdot \alpha_x + \frac{\Delta t}{\Delta y}\cdot \alpha_y \le 1$.
\end{example}

\begin{example}
	An example of numerical hamiltonian $h^A$ satisfying the assumptions required is the \emph{Lax-Wendroff} hamiltonian
	\begin{align}
	\label{LW_2D}
	h^A(x,y,D_x^\pm u,D_y^\pm u)&:=H(x,y,D_x u,D_y u)-\nonumber\\
	&\frac{\Delta t}{2}\left[H_p(x,y,D_x u,D_y u)\left( H_p(x,y,D_x u,D_y u)D^2_{x}u+H_x(x,y,D_x u,D_y u)\right)+\right.\nonumber\\
	&+H_q(x,y,D_x u,D_y u)\left(H_q(x,y,D_x u,D_y u) D^2_{y} u+H_y(x,y,D_x u,D_y u)\right)+\nonumber\\
	&\left.\quad+2H_p(x,y,D_x u,D_y u)H_q(x,y,D_x u,D_y u) D^2_{xy}u\right],
	\end{align}
	where $D^\pm_x u$, $D_x u$, $D^2_{x} u$ are, respectively, the usual one-sided and centered one-dimensional finite difference approximations of the first and second derivative in the $x$-direction (analogously for the $y$-direction), whereas for the mixed derivative we use
	\begin{equation}
	D^2_{xy} u_{i,j} := \frac{u_{i+1,j+1}-u_{i-1,j+1}-u_{i+1,j-1}+u_{i-1,j-1}}{4\Delta x\Delta y}.
	\end{equation}
	Note that the derivatives of $H$ can be computed either analytically or by some second order numerical approximation. In particular, to compute the derivative $H_x$, we can simply use
	\begin{equation}
	(H_x)_{i,j} := \frac{H(x_{j+1},y_i,D_x u_{i,j},D_y u_{i,j})-H(x_{j-1},y_i,D_x u_{i,j},D_y u_{i,j})}{2\Delta x},
	\end{equation}
	and analogously for $H_y$.
\end{example}
For more details on the construction of $S^M$ and $S^A$ and other examples of possible numerical hamiltonians, see \cite{PaolucciPhD,FPT18a}.  

In our approach, in order to couple the two schemes, we need to define three key quantities:
\begin{enumerate}
	\item The \emph{filter function F}, which must satisfy
	\begin{enumerate}
		\item $F(r) \approx r$ for $|r|\leq 1$ so that if $| S^A-S^M|\leq \Delta t \varepsilon^n$ and $\phi_{i,j}^n=1\Rightarrow S^{AF}\approx S^A$,
		\item $F(r) =0$ for $|r|>1$ so that if $| S^A-S^M|> \Delta t \varepsilon^n$ or $\phi^n_{i,j}=0 \Rightarrow S^{AF}= S^M$.
	\end{enumerate}
	Several choices for $F$ are possible, different for regularity properties. In this paper, we will consider the discontinuous filter already used in \cite{BFS16} and defined as follows:
	\begin{equation}
	F(r) := \left\{
	\begin{array}{ll}
	r\qquad&\textrm{ if } |r|\leq 1 \\
	0&\textrm{ otherwise,}
	\end{array}
	\right.
	\end{equation}
	which is clearly discontinuous at $r=-1,1$ and satisfies trivially the two required properties.\\
	\item  If we want the scheme \eqref{AFS_2D} to switch to the high-order scheme when some regularity is detected, we have to choose $\varepsilon^n$ such that
	\begin{equation}
	\label{ep_ineq_2D} 
	\left|\frac{S^A(v^n)_{i,j}-S^M(v^n)_{i,j}}{\varepsilon^n \Delta t}\right|=\left|\frac{h^A(\cdot,\cdot)-h^M(\cdot,\cdot)}{\varepsilon^n}\right| \leq 1,\qquad \textrm{ for }(\Delta t,\Delta x,\Delta y)\to 0,
	\end{equation}
	in the \emph{region of regularity at time $t_n$}, that is  
	\begin{equation}\label{def_region_Rn}
	\mathcal R^n := \left\{(x_j,y_i) : \phi^n_{i,j}=1\right\}.
	\end{equation}
	Proceeding by Taylor expansion for the monotone and the high-order Hamiltonians, by  \eqref{ep_ineq_2D} we arrive to a lower bound for $\varepsilon^n$. 
	The simplest numerical approximation of that lower bound is the following
	\begin{align}\label{eps_2D}
	\varepsilon^n= \max_{(x_j,y_i)\in\mathcal R^n}K & \left|\frac{\Delta t}{2}\left[H_p\left(H_x+H_p D^2_x u^n\right)+H_q\left(H_y+H_q D^2_y u^n\right)+2H_p H_q D^2_{xy} u^n)\right] +\right. \nonumber\\
	&\left. \left(\widetilde h^M_{p^+}-\widetilde h^M_{p^-}\right)+\left(\widetilde h^M_{q^+}-\widetilde h^M_{q^-}\right)\right|,
	\end{align}
	in which we have used the usual notation for the gradient, i.e. $(p,q) := (v_x,v_y)$ and 
	\begin{equation}\label{def:hMp+}
	\widetilde h^M_{p^+} := h^M\left(x,y,D_x u^n,D^+_x u^n, D_y u^n,D_y u^n\right)-h^M\left(x,y,D_x u^n,D^-_x u^n, D_y u^n,D_y u^n\right). 
	\end{equation} 
	The definition of $\widetilde h^M_{p^-}, \widetilde h^M_{q^+}, \widetilde h^M_{q^-}$ follows from \eqref{def:hMp+} in an analogous way. 
	All the derivatives of $H$ are computed at $(x,y,D_x u^n,D_y u^n)$ and the finite difference approximations around the point $(i,j)$, using $K>\frac{1}{2}$. 
	See \cite{PaolucciPhD} for more details.   
	\item For the definition of a function $\phi$, needed to detect the region  $\mathcal R^n$, we require 
	\begin{equation}\label{phi_2D}
	\phi^n_{i,j} :=\left\{
	\begin{array}{ll}
	1\qquad&\textrm{ if the solution } u^n\textrm{ is regular in }I_{i,j},  \\
	0&\textrm{ if }I_{i,j} \textrm{ contains a point of singularity,}
	\end{array}
	\right.
	\end{equation}
	with $I_{i,j} := [x_{j-1},x_{j+1}]\times[y_{i-1},y_{i+1}]$. 
	In order to proceed with the construction, we split the cell $I_{i,j}$ into four subcells, denoted by the superscript `$\vartheta_1 \vartheta_2$', for $\vartheta_1, \vartheta_2=+,-$, according to the shift with respect to the center $(x_j,y_i)$. Then, extending the classical WENO approach proposed in \cite{JP00} to multiple spatial dimensions, we measure the regularity of the solution inside each subcell by computing the \emph{smoothness coefficients} as  rescaled $L^2$ norms of the Lagrange polynomial $P_k^{\vartheta_1,\vartheta_2}(x,y)$ interpolating the values of $u^n$ on the considered stencil, that is
	\begin{align}\label{explicit_beta_2D}
	\beta_k&^{\vartheta_1\vartheta_2}=(-1)^{|\vartheta|}\sum_{{\genfrac{}{}{0pt}{2}{\alpha_1, \alpha_2=0}{|\alpha|\geq 2}}}^2 \int_{\vartheta_1\Delta x}^0\int_{\vartheta_2\Delta y}^{0} \Delta x^{2(\alpha_1-1)}\Delta y^{2(\alpha_2-1)} \left(\partial_x^{\alpha_1}\partial_y^{\alpha_2} P_k^{\vartheta_1 \vartheta_2}(x,y)\right)^2 dx dy\nonumber \\
	&=\frac{1}{\Delta x\Delta y}\left[ {u_{[2,0]}}^2+{u_{[0,2]}}^2+{u_{[1,1]}}^2+\frac{17}{12}\left({u_{[2,1]}}^2+{u_{[1,2]}}^2\right)+\frac{317}{720}{u_{[2,2]}}^2+u_{[2,0]}u_{[2,1]}\right.\nonumber\\
	&\quad\qquad\qquad\left.+u_{[0,2]}u_{[1,2]}-\frac{1}{6}\left(u_{[2,0]}u_{[2,2]}+u_{[0,2]}u_{[2,2]}\right)-\frac{1}{12}\left(u_{[2,1]}u_{[2,2]}+u_{[1,2]}u_{[2,2]}\right)\right] 
	\end{align}
	where $|\vartheta|$ denotes the number of `$-$' in $(\vartheta_1,\vartheta_2)$, for $k=0,1$. Note that we have dropped the dependence on the time step $t^n$ for brevity and we have used the shorter notation $u_{[t,s]}$ to denote the multivariate undivided difference of $u$ of order $t$ in $x$ and $s$ in $y$. The previous formula can be used to obtain all the needed quantities as long as the following \emph{ordered} stencils are used to compute the undivided differences\\
	
	$\bullet$ $\mathcal S_0^{--}=\{x_{j-1},x_j,x_{j+1}\}\times\{y_{i-1},y_i,y_{i+1}\}$, $\mathcal S_1^{--}=\{x_{j},x_{j-1},x_{j-2}\}\times\{y_{i},y_{i-1},y_{i-2}\}$;\\
	$\bullet$ $\mathcal S_0^{+-}=\{x_{j+1},x_{j},x_{j-1}\}\times\{y_{i-1},y_i,y_{i+1}\}$, $\mathcal S_1^{+-}=\{x_{j},x_{j+1},x_{j+2}\}\times\{y_{i},y_{i-1},y_{i-2}\}$;\\
	$\bullet$ $\mathcal S_0^{++}=\{x_{j+1},x_{j},x_{j-1}\}\times\{y_{i+1},y_i,y_{i-1}\}$, $\mathcal S_1^{++}=\{x_{j},x_{j+1},x_{j+2}\}\times\{y_{i},y_{i+1},y_{i+2}\}$;\\
	$\bullet$ $\mathcal S_0^{-+}=\{x_{j-1},x_{j},x_{j+1}\}\times\{y_{i+1},y_i,y_{i-1}\}$, $\mathcal S_1^{-+}=\{x_{j},x_{j-1},x_{j-2}\}\times\{y_{i},y_{i+1},y_{i+2}\}$.\\
	
	Since these coefficients are such that\\
	
	$\bullet$ $\beta_k=O(\Delta^2)$, with $\Delta:=\max\{\Delta x,\Delta y\}$, if the solution is smooth in $\mathcal S_{k}$;\\
	$\bullet$ $\beta_k=O(1)$ if there is a singularity in $\mathcal S_{k}$,\\
	
	according to the usual WENO procedure we weight the obtained information and focus on the `inner' stencil, denoted by the subscript `$0$', by computing
	\begin{equation}
	\omega^{\vartheta_1\vartheta_2}=\frac{\alpha^{\vartheta_1\vartheta_2}_0}{\alpha^{\vartheta_1\vartheta_2}_0+\alpha^{\vartheta_1\vartheta_2}_1},
	\end{equation}
	where $\alpha^{\vartheta_1\vartheta_2}_k=\frac{1}{(\beta^{\vartheta_1\vartheta_2}_k+\sigma_\Delta)^2}$,  
	with $\sigma_\Delta=\Delta x^2+\Delta y^2$, 
	which represents the measure of smoothness of the solution in each subcell. Once we have computed the four indicators, we couple the information by defining
	\begin{equation}
	\label{omega_2D}
	\omega=\min\{\omega^{--},\omega^{+-},\omega^{-+},\omega^{++}\},
	\end{equation}
	which, as can be shown by exploiting the properties of the coefficients $\beta_k$, is such that
	\begin{equation}
	\omega_{i,j}=
	\left\{
	\begin{array}{ll}
	O(\Delta ^4)\qquad &\textrm{ if }\rho \in I_{i,j} \\
	\frac{1}{2}+O(\Delta)& \textrm{ otherwise,}
	\end{array}
	\right.
	\end{equation}
	where $\rho$ is a point of discontinuity in the gradient. 	
	At this point, in order to reduce the amplitude of the oscillations around the optimal value $\frac{1}{2}$ in regular regions (the $O(\Delta)$ term), we use the \emph{mapping} first introduced in \cite{HAP05} to propose a modification of the original WENO procedure, called \emph{M-WENO}, 
	that is
	\begin{equation}
	\label{map2}
	\omega^*=g(\omega)=4\omega\left(\frac{3}{4}-\frac{3}{2}\omega+\omega^2\right),
	\end{equation}
	which, using Taylor expansion around $\frac{1}{2}$, directly gives $\omega^*=\frac{1}{2}+O(\Delta^3)$ when the solution is regular in $I_{i,j}$.
	Finally, in order to define our function $\phi$, it is enough to take
	\begin{equation}
	\label{phi_disc}
	\phi(\omega^*)=\chi_{\{\omega^*\geq M\}},
	\end{equation}
	with $M<\frac{1}{2}$ (e.g. $M=0.1$ in the tests reported in Sect. \ref{sec:tests}), a number that can also depend on $\Delta x$ and $\Delta y$ . 
	For more details and other possible constructions for the definition of $\phi$, we refer the interested reader to \cite{PaolucciPhD,FPT18a}. 
\end{enumerate}
\begin{remark}
	Choosing $\varepsilon^n\equiv \varepsilon\Delta x$, with $\varepsilon>0$ and $\phi_{i,j}^n\equiv 1$, we get the Filtered Scheme of \cite{BFS16}, so here we are generalizing that approach to exploit more carefully the local regularity of the solution at every time $t^n$ and cell $I_{i,j}$.
\end{remark}

\section{Numerical implementation of the modified LS method}\label{sec:num_impl}
Before illustrating the numerical tests, let us first give some comments on the numerical schemes composing the AF scheme adopted for the tests in Sect. \ref{sec:tests}. The main issue concerning the \emph{local Lax-Friedrichs} and the \emph{Lax-Wendroff schemes} defined by (\ref{local_lax_fried_2D}) and (\ref{LW_2D}), respectively, is the need to compute the one-directional velocities $H_p$ and $H_q$ which depend also on $\frac{\partial c}{\partial \xi}$ and $\frac{\partial c}{\partial \zeta}$, as visible in (\ref{def_Hp}). Moreover, in order to implement the local Lax-Friedrichs scheme we should be able to compute the maximum of $|H_p|$ (resp. $|H_q|$) uniformly with respect to $p$ (resp. $q$), 
which is a very intricate matter due to the (possible) low regularity of $\widetilde c$. 
In fact, if we focus on the usual behavior of $c(x,y)$ in the proximity of a relevant edge, we can expect the derivatives $\frac{\partial c}{\partial \xi}$ and $\frac{\partial c}{\partial \zeta}$ to be really big. This is not surprising since the front decelerates rapidly in the neighborhood of an edge. 
In addition, in order to solve the full model (\ref{eq:HJ_eikonal})-(\ref{eq_new_vel}), we should take into account also the remaining dependence of $H(x,y,v,\cdot,\cdot)$ when deriving the second-order Lax-Wendroff scheme and, clearly, the formula to compute the threshold $\varepsilon^n$. 
Finally, concerning the Courant-Friedrichs-Lewy (CFL) condition, it is necessary to compute $\max\{|H_p|,|H_q|\}$ with the full formula (\ref{def_Hp}). Consequently, 
$\lambda$ could be excessively small due to the low regularity of $\widetilde c$.  
In this latter case, we would clearly need an adaptive mesh refinement technique to reduce the computational cost.

In order to avoid most of these complications in the numerical implementation, we choose to approximate the solution of the simplified problem (\ref{eq_simplified}), adjusting the velocity $\widetilde c$ according to (\ref{eq_new_vel}) at each time step. 
Using the simplified problem \eqref{eq_simplified}, we can use the simple relation
\begin{equation}\label{simple_relation}
\max_{p}\max_{q} |H_p(\cdot,p,q)|=\max_{p} |H_p(\cdot,p,0)|,
\end{equation}
avoiding to take the maximums over all the possible values of $p$ and $q$, which is instead required for the resolution of the full problem (\ref{eq:HJ_eikonal})-(\ref{eq_new_vel}) with anisotropic velocity $\widetilde c$. 
Analogous comment holds for $H_q$. 
Lastly, from the numerical point of view, the use of this  simplification brings another fundamental consequence: when we apply the numerical schemes to solve (\ref{eq_simplified}), we are considering, formally, a problem with bounded velocities $\max\{|H_p|,|H_q|\}\leq 1$. This implies that we can choose the following CFL condition:
\begin{equation}\label{CFL_segm}
\lambda:=\max\left\{\frac{\Delta t}{\Delta x},\frac{\Delta t}{\Delta y}\right\}\leq \frac{1}{2}\max\{|H_p|^{-1},|H_q|^{-1}\},
\end{equation}
using the relation (\ref{simple_relation}), which is a less restrictive condition with respect to the original one coming from the full problem (\ref{eq:HJ_eikonal})-(\ref{eq_new_vel}). 

In the following, we will use the same notations introduced in Sect. \ref{sec:AFS}, except for the number of time steps $N_T$, which will be replaced by the total number of iterations $N_i$ used by the scheme, since now we are looking for an asymptotic solution (in some stationary sense). The maximum number of iterations, which is fixed at the beginning of the procedure, will be denoted by $N_{\max}$.

Let us give some details on the precise numerical implementation, commenting the main procedures involved in the (sketched) \emph{Algorithm} \ref{pseudo_code}.\\ 
\begin{algorithm}[h!]
	\caption{Segmentation via the LS Method\label{pseudo_code}}
	\begin{algorithmic}
	\Require{$\mu$, $K_{reg}$, $tol$, $N_{\max}$, $I$, $u_0$}
\State $E^0=1$, $n=0$
\State regularize the matrix $I$ (apply the Gaussian filter)
\State compute the velocity matrix $c$ using (\ref{vel_c1}) or \eqref{vel_c2} 
\State store the position of the front in the matrix $F^0$

\While {$(E^n\ge tol)$ and $(n<N_{\max})$}
\State Step 1: compute the modified velocity matrix $\widetilde c^n$ using (\ref{eq_new_vel})
\State Step 2: update the solution $u^{n}\to u^{n+1}$
\State $n=n+1$
\State Step 3: store the front $F^n$
\State  \qquad\quad \, compute the error $E^n$ 
\EndWhile

\State $N_i=n$
\Ensure{$N_i$, $u^{N_i}$.}
\end{algorithmic}
\end{algorithm}

Let us set the parameters of the simulation, which are the power $\mu$ in (\ref{vel_c1}), 
the number of iterations $K_{reg}$ of the heat equation for the Gaussian filter, the tolerance $tol>0$ of the stopping criterion, the amplitude of the pixels $(\Delta x,\Delta y)$ and, subsequently, the time step $\Delta t$ according to the CFL condition (\ref{CFL_segm}).

Then, at each iteration $n=0,\dots,N_i$, which has to be interpreted in the sense ``until convergence'' (note that $N_i$ is not known a priori, but depends on the stopping criterion described in Step 3 and can be equal to $N_{max}$ in case of not convergence of the scheme), we repeat the following steps.

\noindent
\emph{Step 1.} 
For $i=0,\cdots N_y, j=0,\cdots, N_x$, with $(N_x+1)\times (N_y+1)$ the size of the input image, we precompute the matrix $\widetilde c(x_j,y_i,u_{i,j},D_x u_{i,j},D_y u_{i,j})$ at the beginning of each iteration using central finite difference approximations for the first order derivatives $D_x u_{i,j}$ and $D_y u_{i,j}$. 
Note that the quantities only depend on $(i,j)$ also through $u$. Clearly, this method is valid only as long as the representation function $u$ remains smooth at all the level sets, and should be justified in the case of singular edges (although we will not pursue this precise matter). Moreover, in general the point
\begin{equation}
\left(x_{j_u},y_{i_u}\right):=\left(x_j-d(u_{i,j})\frac{D_x u_{i,j}}{\sqrt{{(D_x u_{i,j})}^2+{(D_y u_{i,j})}^2}},y_i-d(u_{i,j})\frac{D_x u_{i,j}}{\sqrt{{(D_x u_{i,j})}^2+{(D_y u_{i,j})}^2}}\right)
\end{equation}
is not a point of the grid $(x_j,y_i)$. \\
To reconstruct the correct value (or at least a reasonable approximation) there exist different possible implementations. For example, a simple \emph{bilinear reconstruction} from the neighboring values 
\begin{equation}
\mathcal N_u:=\left\{\left(x_{\lfloor j_u\rfloor},y_{\lfloor i_u \rfloor}\right), \left(x_{\lceil j_u\rceil},y_{\lceil i_u \rceil}\right),\left(x_{\lceil j_u\rceil},y_{\lfloor i_u \rfloor}\right),\left(x_{\lfloor j_u\rfloor},y_{\lceil i_u \rceil}\right)\right\},
\end{equation}
where we have used the notation
\begin{equation}
\lceil j_u\rceil:=j-\left\lceil\frac{ x_{j_u}-x_j}{\Delta x}\right\rceil\quad \textrm{ and } \quad \lfloor i_u\rfloor:=i-\left\lfloor \frac{y_{i_u}-y_i}{\Delta y}\right\rfloor,
\end{equation}
with the other cases following an analogous definition.
Another possibility, which we have used in the numerical tests since it seems to give nicer results in terms of the shape of the approximate representation $u$, consists 
in taking as $\left(x_{j_u},y_{i_u}\right)$ the point such that 
\begin{equation}
|u_{i_u,j_u}|:=\min_{(x_j,y_i)\in\mathcal N_u} |u_{i,j}|.
\end{equation}
Note that this construction is well defined only if $|\nabla u_{i,j}|\not=0$. Therefore, we define the updated velocity matrix as 
\begin{equation}
\widetilde c^n_{i,j} :=\left\{
\begin{array}{ll}
c_{i_u,j_u}\qquad&\textrm{ if } |\nabla u^n_{i,j}|\not=0,\\
c_{i,j}&\textrm{ otherwise,}
\end{array}
\right.
\end{equation}
and we use $\widetilde c^n$ as an isotropic velocity in the next step.

\noindent
\emph{Step 2.} We approximate the problem \eqref{eq_simplified} 
using the AF scheme (\ref{AFS_2D}), with the local Lax-Friedrichs scheme (\ref{local_lax_fried_2D}) as $S^M$ and the Lax-Wendroff scheme (\ref{LW_2D}) as $S^A$. 
From now on, we will refer to this AF implemented scheme with the acronym AF-LW.  \\
We add homogeneous Neumann boundary conditions to the problem \eqref{eq_simplified} in all our experiments  in order to not alter the average intensity of the image.

\noindent
\emph{Step 3.} In this last step, we describe how approximating the front $\Gamma_t$. Since $\Gamma_t$  is a curve, it is composed by points that are not all grid points belonging to our mesh. Hence, in order to approximate the position of the front at each time step $t$, we consider a neighborhood $\theta_\delta$ of the front $\Gamma_t$ of radius $\delta=\max\{\Delta x,\Delta y\}$. In this way, stopping the evolution as soon as the front ceases to move will be equivalent to require that the neighborhood $\theta_\delta$ ceases to move. 
In order to apply that procedure, at each iteration $n$ we store the values of the points $(x_j,y_i)$ such that $u_{i,j}^n$ changes sign in a matrix $F_{i,j}^n$  (we use the closest points on the grid, that are $(i,j\pm1)$ and $(i\pm1,j)$), 
and set $F^n_{i,j}=0$ otherwise. 
In this way we automatically store the disposition of the front with an error of order $\delta=\max\{\Delta x,\Delta y\}$, i.e. we approximate the position of the front, as desired. \\
We will continue to do that at each iteration until the matrix $F$ at two consecutive iterations will be ``close enough". For this reason, we consider two stopping rules: 
\begin{equation}\label{eq:err_inf}
E_\infty:=||u^{n+1}-u^n||_{L^\infty(\theta_\delta)}=\max_{i,j}|F^n_{i,j}-F^{n-1}_{i,j}|<\tau,
\end{equation}
where $\tau>0$ is the prescribed tolerance a priori chosen, 
\begin{equation}\label{eq:err_uno}
E_1:=||u^{n+1}-u^n||_{L^1(\theta_\delta)}=\Delta x\Delta y \sum_{i,j}|F^n_{i,j}-F^{n-1}_{i,j}|<\tau, 
\end{equation}
where now a dependence on the discretization parameters appears. \\ 
In our implementation we use one of the two \emph{stopping rules} above introduced combined with a condition on the number of allowed iterations (i.e. $n<N_{max}$). 

\section{Numerical simulations}\label{sec:tests}
In this section we present a series of numerical experiments on both synthetic and real images, comparing the results obtained by the AF scheme with those obtained by the simple monotone scheme and a high-order scheme which uses the Total Variation Diminishing (TVD) Runge-Kutta (RK) of third order in time and the WENO scheme of second/third order in space. 
More in details, for the WENO scheme we use the same efficient implementation suggested in \cite{JP00}, Remark 1 on page 2130, which we extended to the 2D case, adding the improvement presented in \cite{ABM10}, which consists in choosing $\sigma=\Delta x^2$ instead of $\sigma=10^{-8}$ that is the value used in the original paper \cite{JP00} by Jiang and Peng. 
Note that we used this scheme instead of the third to fifth order scheme used in the numerical tests of \cite{JP00} for comparison reasons, otherwise it would be not comparable by order. 
The first aim is to show the possible improvements of the modified model with respect to the classical formulation. In fact, after extensive numerical simulations, we noticed that the classical model is not well defined when using high-order schemes, since they can produce heavy oscillations 
as soon as the Lipschitz constant of the representation function becomes too big. 
This effect causes the stopping rule to be practically ineffective (independently of the norm used) in most cases particularly when using the AF scheme or the WENO scheme, whereas the simple monotone scheme seems to give always stable results. 
Note that, when the singularity develops, the representation function becomes more and more vertical as time flows. 

The numerical tests illustrated in this section will compare the results also varying the initial datum or varying the norm for the stopping rule defined in \eqref{eq:err_inf} and \eqref{eq:err_uno}. 
Moreover, for synthetic images, we also vary the space steps $\Delta x$ and $\Delta y$, that we consider equal to each other ($\Delta x=\Delta y$), and, therefore, the total number of pixels. 

Now, let us specify the initial condition used in each case. When the velocity is defined by the classical model, in the expansion case (\emph{Case a}) we use the paraboloid
\begin{equation}\label{parab_u0}
u_0(x,y)=\min\left\{x^2+y^2-r^2,\frac{1}{2}r^2\right\},
\end{equation}
where $r>0$ is the radius of the initial circle and $\frac{1}{2}r^2$ a value chosen in order to cut the surface from above (therefore we have a flat surface at the numerical boundary),whereas in the shrinking case (\emph{Case b}) we use the  truncated pyramid  (with a square or rectangular base depending on the image frame), that is 
\begin{equation}\label{truncated_pyramid_u0}
u_0(x,y)=\min\{2(x-b_x),2(a_x-x),2(y-b_y),2(a_y-y),-0.2\},
\end{equation}
where $[a_x,b_x]\times[a_y,b_y]$ is the frame of the image, $-0.2$ is the value at which we truncate the pyramid and $2$ is the steepness of the faces of the surface 
(the smaller the value, the less steep the faces are). 
By this choice we use a slightly more regular front with respect to the discontinuous representation that simply changes value crossing the frame of the image, still being able to keep the whole surface outside the region occupied by the objects to segment.

For the modified velocity $\widetilde c$, in the expansion case (\emph{Case a}) we consider two different initial data: the paraboloid as in (\ref{parab_u0}) (\emph{Datum 1}) or the following signed distance function (\emph{Datum 2})
\begin{equation}\label{distance_u0_exp_case}
u_0(x,y)=dist\{(x,y),\Gamma^1_0\},
\end{equation}
where $\Gamma^1_0$ is 
the usual circle centered in $(0,0)$ with radius $r=0.5$ unless otherwise stated.  
Instead, for the shrinking case (\emph{Case b}) we only use the signed distance function 
\begin{equation}\label{distance_u0_shrink_case}
u_0(x,y)=dist\{(x,y),\Gamma^2_0\},
\end{equation}
with $\Gamma^2_0$ representing the frame of the image. 
We summarize the different initial conditions in Tab. \ref{table:summary_u0}.
\begin{table}[h!] 
	\caption{Summary of the initial conditions considered in the numerical tests. \label{table:summary_u0}}
	\centering
	\begin{tabular}{c|c|c  }
		& $c$ & $\widetilde c$ \\
		\hline
		Expansion case  & $u_0$ as Paraboloid \eqref{parab_u0} &  Datum 1: $u_0$ as Paraboloid \eqref{parab_u0} \\
		(Case a)		&						& Datum 2: $u_0$ as  Distance from $\Gamma^1_0$ \eqref{distance_u0_exp_case} \\
		\hline
		Shrinking case & $u_0$ as Truncated pyramid &      \\
		(Case b)		&	(or tent) \eqref{truncated_pyramid_u0} 	& $u_0$ as Distance from $\Gamma^2_0$ \eqref{distance_u0_shrink_case} \\
		\hline
	\end{tabular}
\end{table}

In order to give a quantitative evaluation of the performances in addition to the qualitative analysis, in the tables we compare the results in terms of \emph{number of iterations} $N_i$ and \emph{relative error in pixels}, defined as 
\begin{equation}\label{eq:err-rel-formula}
P\textrm{-}Err_{rel}=\frac{|P_{ex}-P_{a}|}{P_{ex}},
\end{equation}
where $P_{ex}$ and $P_a$ are the number of pixels inside the exact and approximated boundaries of the object(s), respectively. Note that we can compute the ``exact'' object  
only if the background is really smooth (in the synthetic cases it is always uniform), because we usually use a comparison with a ``threshold'' for the values of $I(x,y)$ in order to select the regions occupied by the object (exact object).  
Whereas, for the approximated object we will count the pixels for which $u_{i,j}\leq 0$. 
Moreover, we measure the error also with a closely related quantity, that is
\begin{equation}\label{eq:err1-formula}
P\textrm{-}Err_1=|P_{ex}-P_{a}|\Delta x \Delta y,
\end{equation}
in order to show some dependence on the discretization parameters. 

\noindent If the schemes do not converge in the fixed maximum number of allowed iterations $N_{\max}$, we will put a ``$-$'' inside the tables, in place of $N_i$. For all our tests, we will set  $N_{\max}=2000$. Moreover, in case the front does not stop correctly on the boundary of the object, thus giving an unstable and unusable result, we will put an ``X'' in correspondence of the errors column.
For each test, we specify all the values of the parameters involved ($\mu$, $K_{reg}$, $tol$ and $\Delta x=\Delta y$ or $\# Nodes$), the norm used in the stopping rule and the chosen velocity function. 
For all the numerical tests presented in this paper, we use CFL number $\lambda=\max\left\{\frac{\Delta t}{\Delta x},\frac{\Delta t}{\Delta y}\right\}=\frac{1}{2}$, $K=1$ in the formula \eqref{eps_2D} for the computation of $\varepsilon$, and the velocity function $c_1$ defined in \eqref{vel_c1} for the classical model, referring to it simply as the classical $c$. 
All the numerical tests have been implemented in language C++, with plots and computation of the errors in MATLAB. 
The computer used for the simulations is a Notebook Asus F556U Intel Core i7-6500U with speed of 2.59 GHz and 12 GB of RAM.	

\subsection{Synthetic tests}\label{sec:ST}
Let us begin by a simple synthetic example. The main goal here is to compare the behavior of the three schemes with respect to the use of the classical velocity $c$ and the modified one $\widetilde c$, and the performances varying the number of grid nodes. 

\paragraph{Test 1. Rhombus}
For this first test, we perform the simulations only in the case of an expansion, since the results do not vary much in the shrinking case.  

\begin{figure}[h!] 
	\centering
	\includegraphics[width=0.32\textwidth]{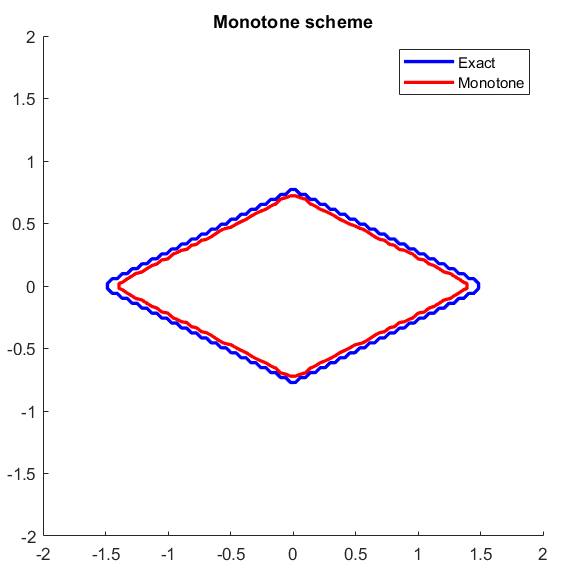}
	\includegraphics[width=0.32\textwidth]{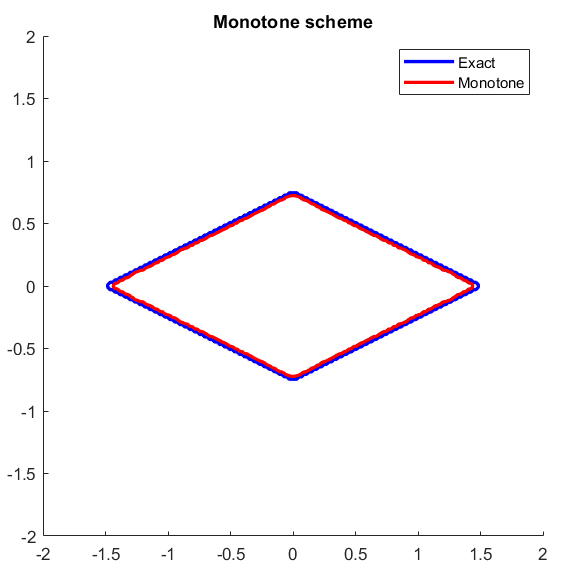}
	\includegraphics[width=0.32\textwidth]{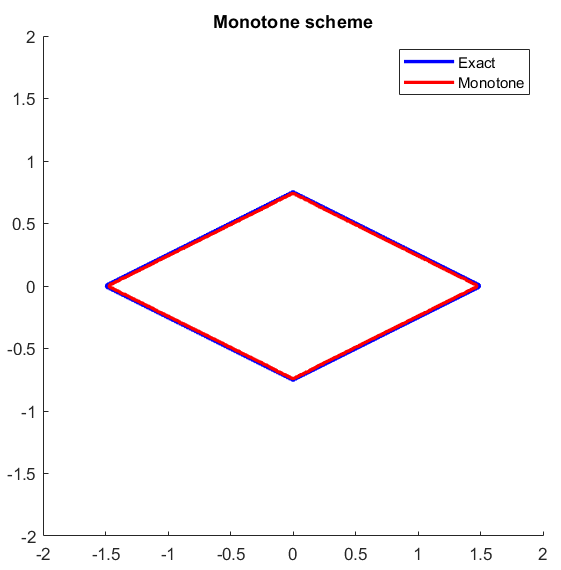}\\
	\includegraphics[width=0.32\textwidth]{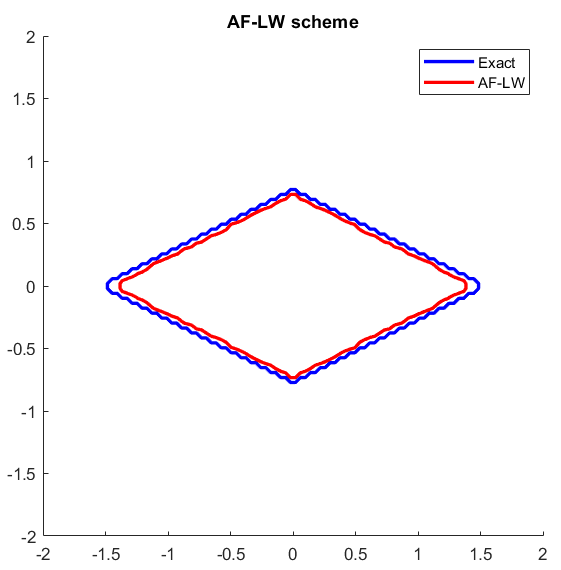}
	\includegraphics[width=0.32\textwidth]{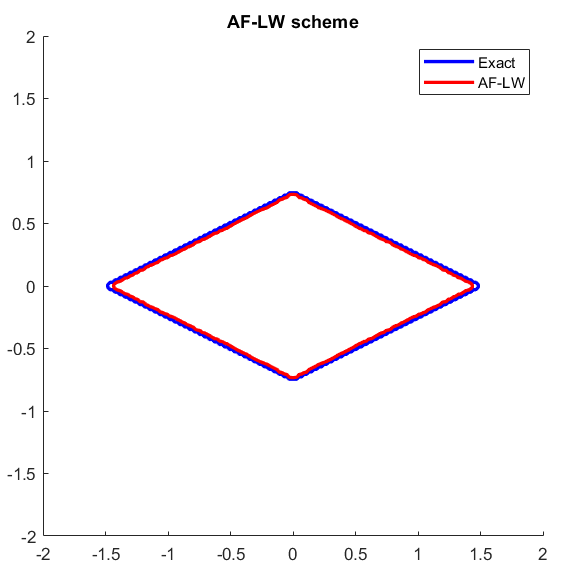}
	\includegraphics[width=0.32\textwidth]{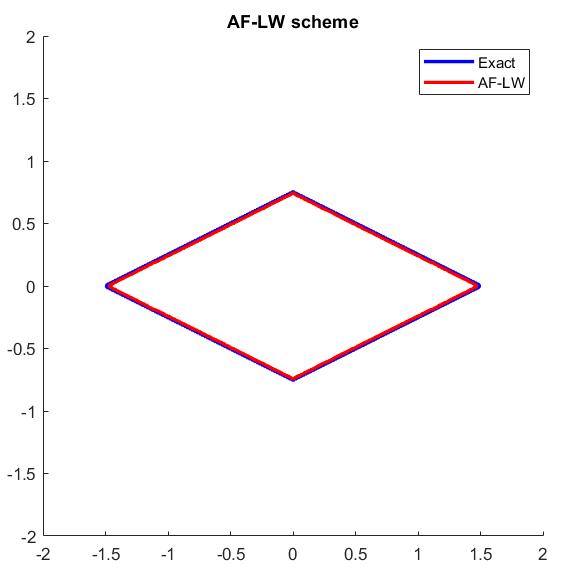}\\
	\includegraphics[width=0.32\textwidth]{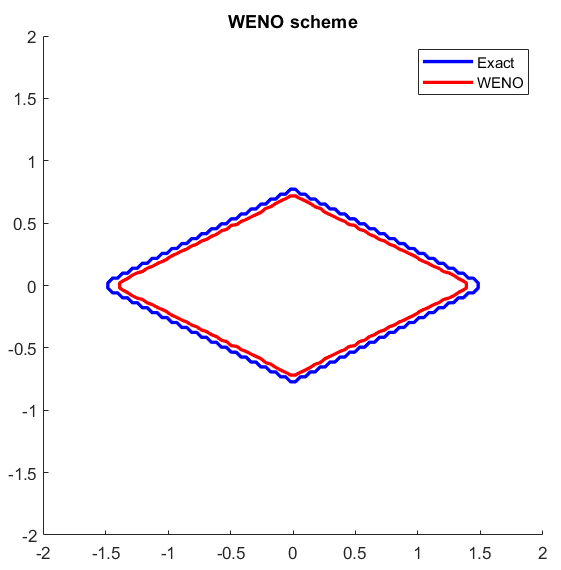}
	\includegraphics[width=0.32\textwidth]{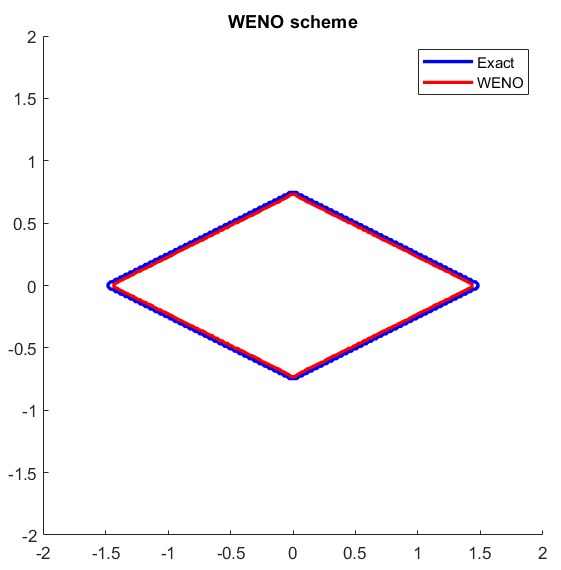}
	\includegraphics[width=0.32\textwidth]{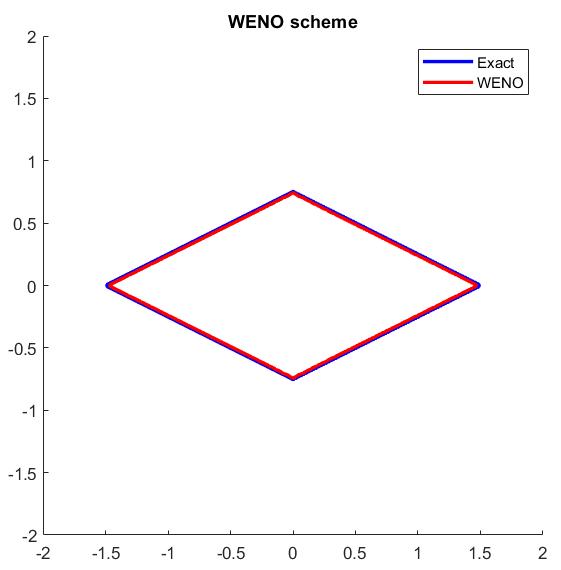}
	\caption{\small{Test 1a with Datum 1. Plots of the final front obtained by the Monotone scheme (top), the AF-LW scheme (middle) and the WENO scheme (bottom) with velocity $\widetilde c$ and the parameters reported in Tab. \ref{table:rhombus_tildec_1_Nx}.}
		\label{fig:test1a}}
\end{figure}

\begin{figure}[h!] 
	\centering
	\includegraphics[width=0.32\textwidth]{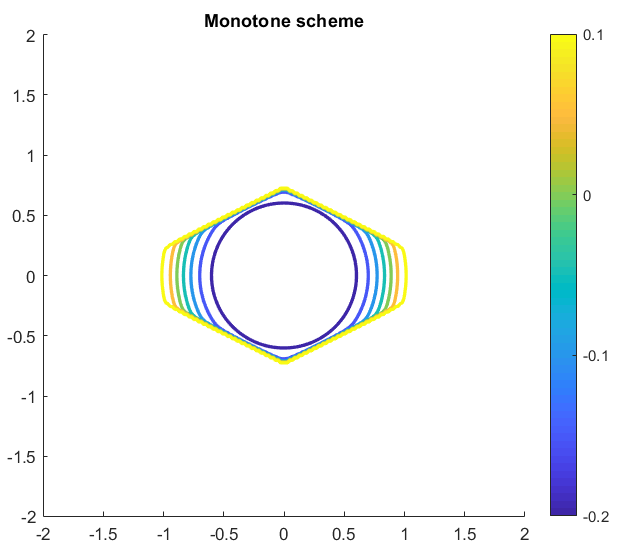}
	\includegraphics[width=0.32\textwidth]{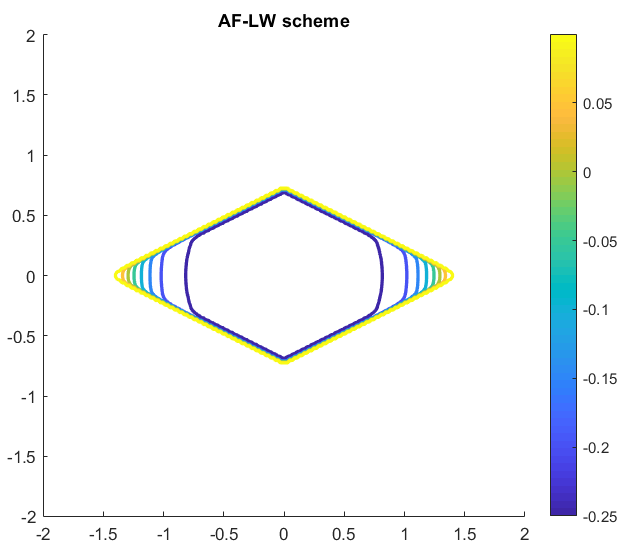}
	\includegraphics[width=0.32\textwidth]{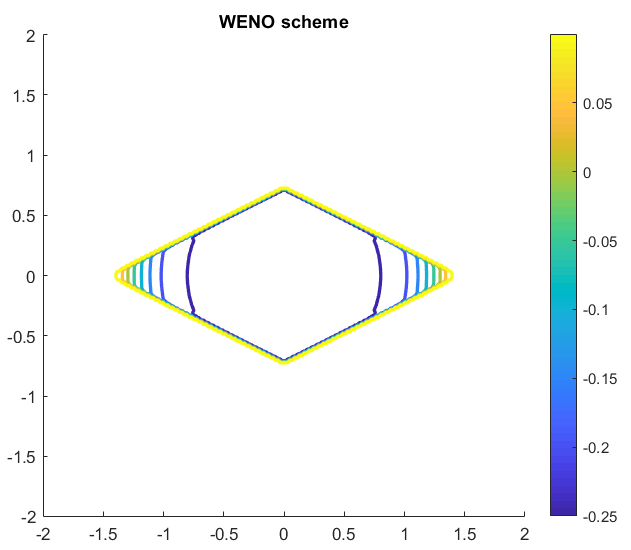}\\
	\includegraphics[width=0.32\textwidth]{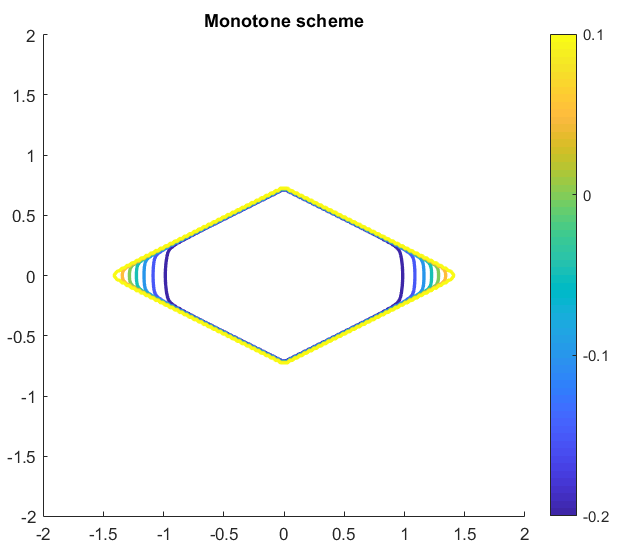}
	\includegraphics[width=0.32\textwidth]{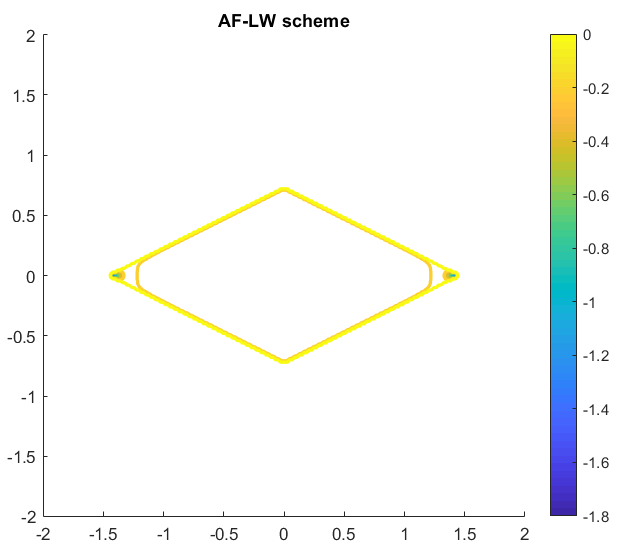}
	\includegraphics[width=0.32\textwidth]{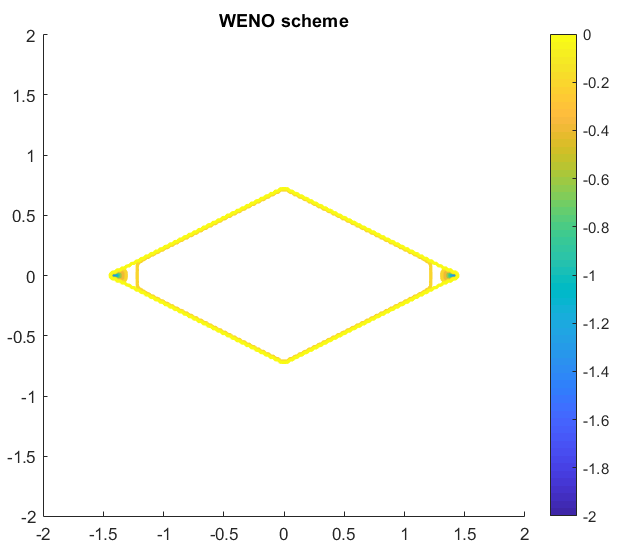}\\
	\includegraphics[width=0.32\textwidth]{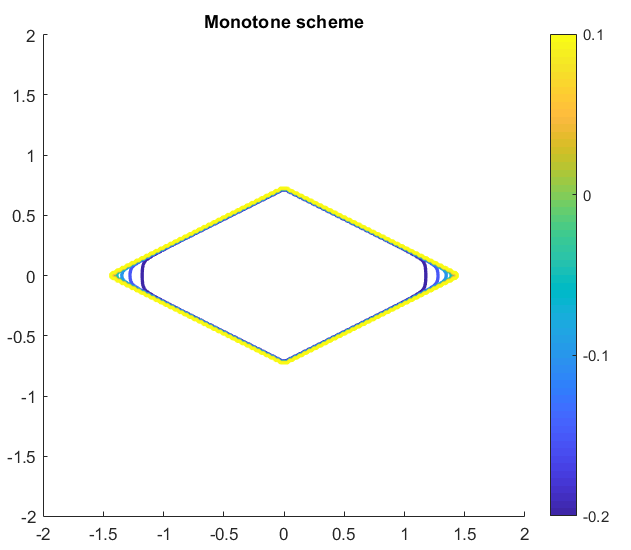}
	\includegraphics[width=0.32\textwidth]{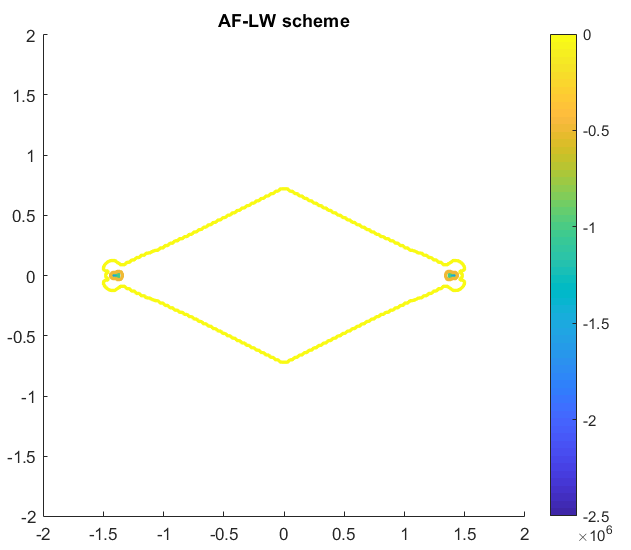}
	\includegraphics[width=0.32\textwidth]{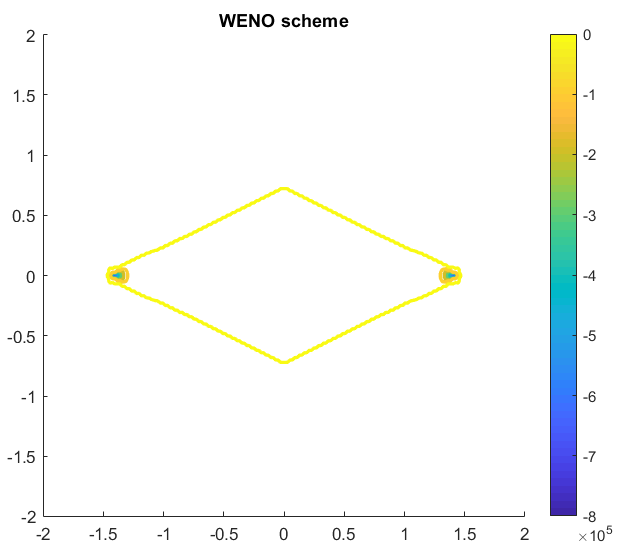}\\
	\includegraphics[width=0.32\textwidth]{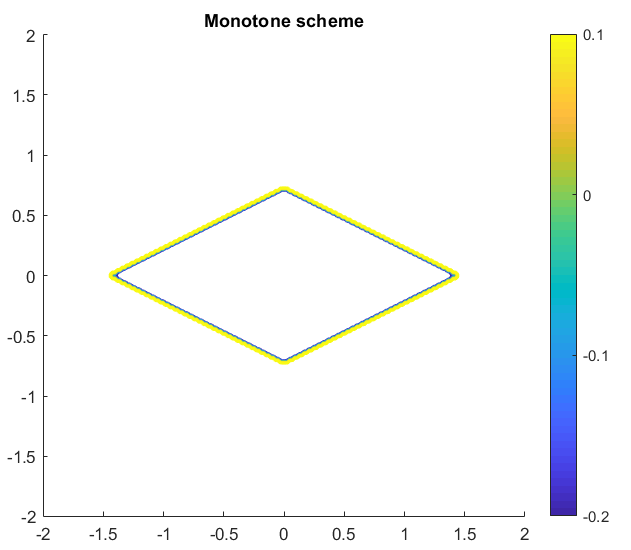}
	\includegraphics[width=0.32\textwidth]{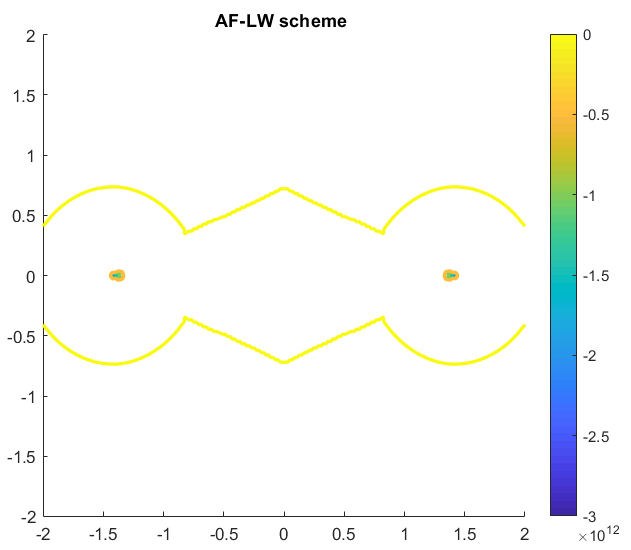}
	\includegraphics[width=0.32\textwidth]{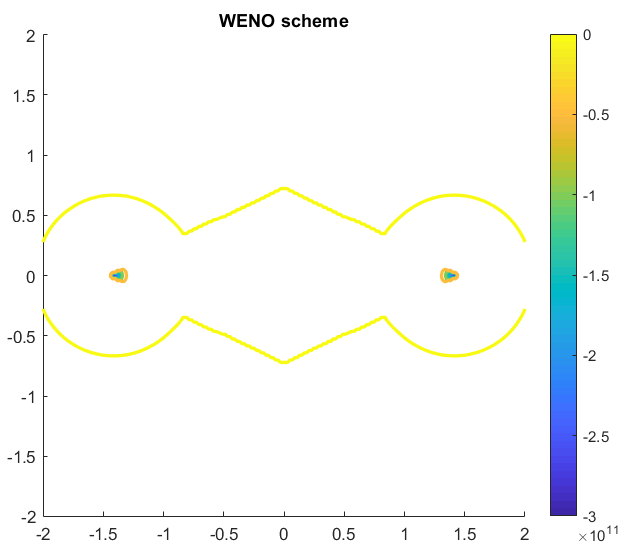}
	\caption{\small{Test 1a with Datum 1. Contour plots of the representation function obtained by the Monotone scheme (left) at $N_i=40$, $80$, $120$ and at final time, and by the AF-LW scheme (middle) and the WENO scheme (right) at $N_i=80$, $100$, $160$, $220$, using velocity $c$ with $\# Nodes=202$.}
		\label{fig:test1a_c_repres}}
\end{figure}

\noindent The rhombus considered is defined by the equation
\begin{equation}
\frac{|x|}{2}+|y|=\frac{3}{4},\qquad(x,y)\in[-2,2]^2,
\end{equation}
that produces a final front visible in Fig. \ref{fig:test1a} for each scheme using the velocity $\widetilde{c}$. As clearly visible, this rhombus presents some heavily marked corners, which causes some serious troubles when using the filtered scheme or the WENO scheme for the solution of the classical model. In fact, looking at Tab. \ref{table:rhombus_c_1_Nx}, we can see that the two high-order schemes achieve convergence, in the sense of the iterative stopping rule \eqref{eq:err_inf}, but after having lost track of the boundary (the front overcomes the edge of the object and then keeps expanding). This happens also varying the number of grid nodes (from $102$ to $402$). 
On the contrary, the monotone scheme converges, the two errors reported decrease by refining, still using  the same parameters $\mu,K_{reg}$ and same tolerance $tol$. The difficulties of the high-order schemes with the classical velocity $c$ are cleary visible looking at Fig. \ref{fig:test1a_c_repres}, in which the contour plots of the representation function obtained by the three schemes are visible: the evolution of the level-sets for the monotone scheme expands and finally coincides at the final time (last row on the left) with the boundary of the object. Instead, for the two high-order schemes we can see different times in order to show well when the oscillations on the left and right edges begin to increase causing the loss of the rhombus boundary in the approximation. 
\begin{figure}[h!] 
	\centering
	\includegraphics[width=0.32\textwidth]{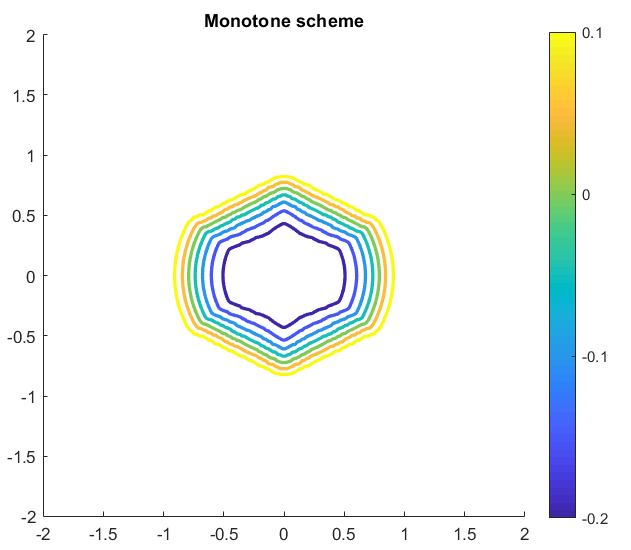}
	\includegraphics[width=0.32\textwidth]{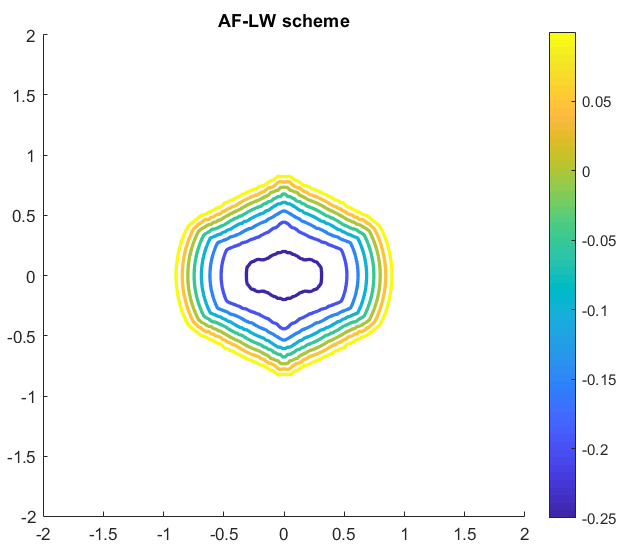}
	\includegraphics[width=0.32\textwidth]{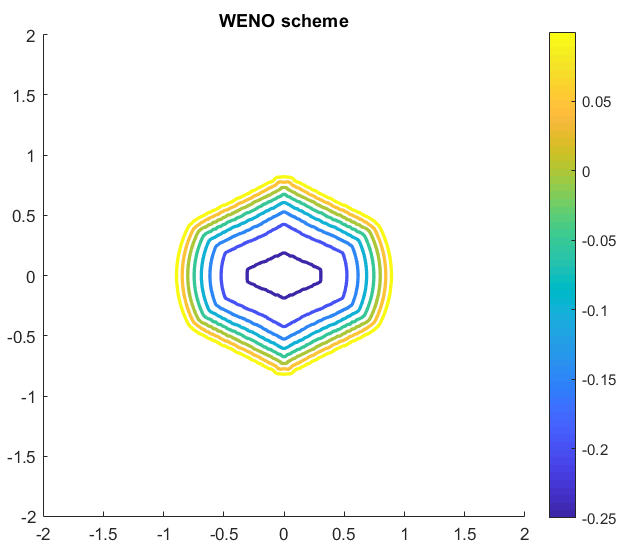}\\
	\includegraphics[width=0.32\textwidth]{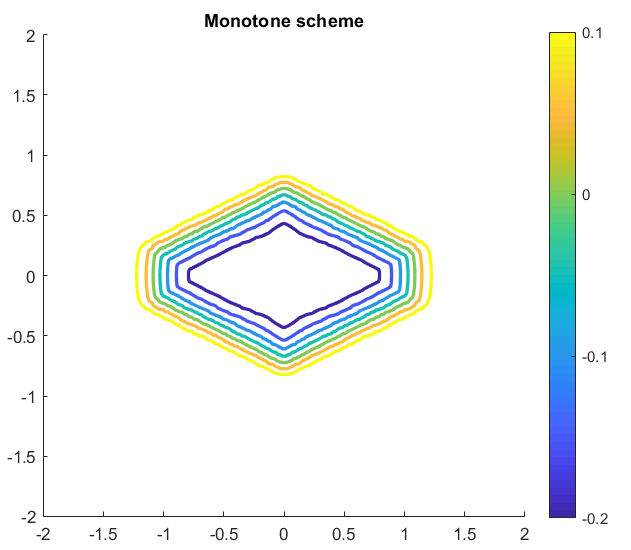}
	\includegraphics[width=0.32\textwidth]{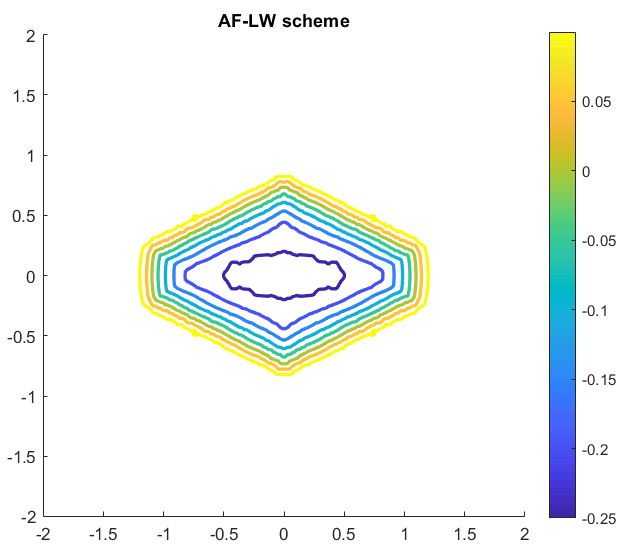}
	\includegraphics[width=0.32\textwidth]{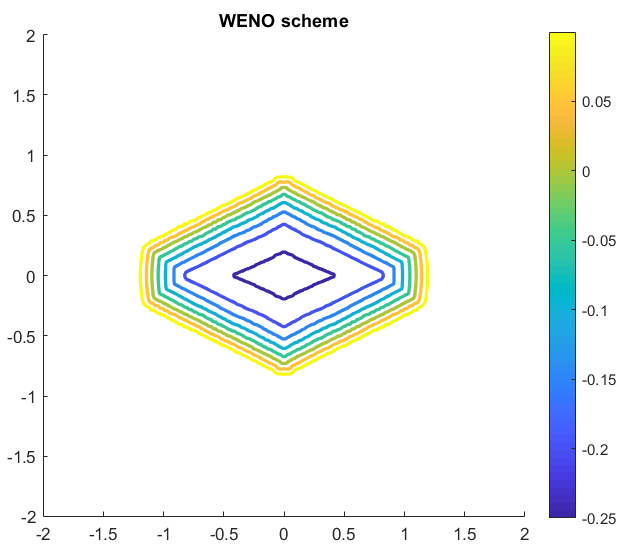}\\
	\includegraphics[width=0.32\textwidth]{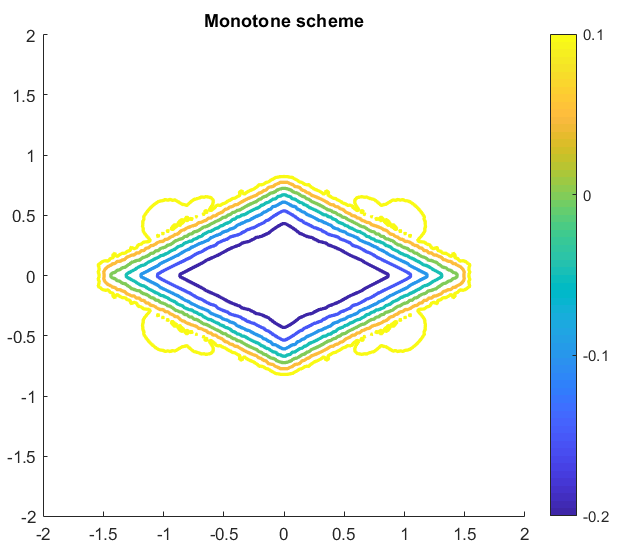}
	\includegraphics[width=0.32\textwidth]{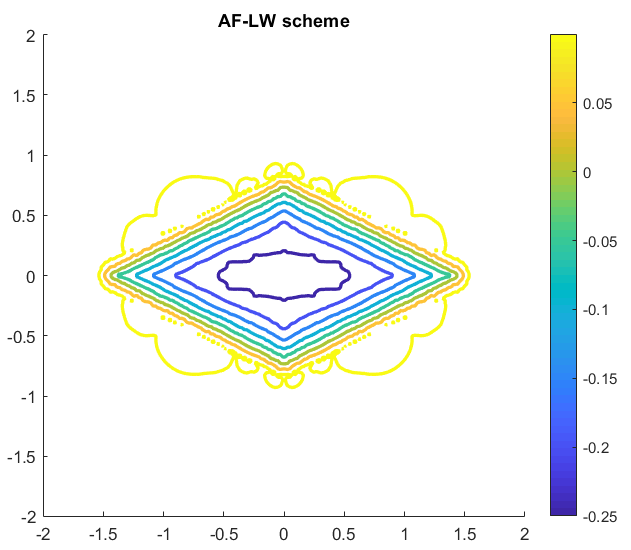}
	\includegraphics[width=0.32\textwidth]{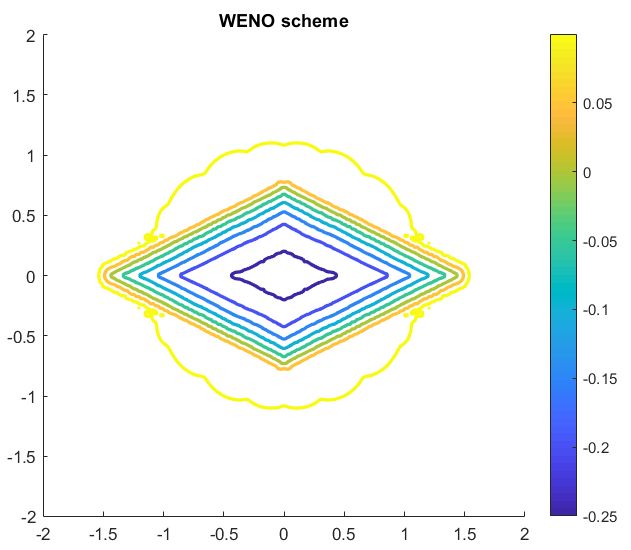}
	\caption{\small{Test 1a with Datum 1. Contour plots of the representation function obtained by the Monotone scheme (left), the AF-LW scheme (middle) and the WENO scheme (right) at $N_i=30$, $60$ and at final time, using velocity $\widetilde{c}$ with $\# Nodes=202$.}
		\label{fig:test1a_ctilde}}
\end{figure}
For the modified velocity $\widetilde c$, the behavior of the schemes is different, see the contour plots in Fig. \ref{fig:test1a_ctilde}. Note that the yellow level-set on the last row, particularly for the WENO scheme, is not the 0-level set. 
The effects of the modified velocity on the evolution is clear comparing Figs. \ref{fig:test1a_c_repres} and \ref{fig:test1a_ctilde}. In the first case all the level sets expand towards the boundary of the rhombus, eventually collapsing onto each other. This gives rise  to a discontinuity around the front, which causes instabilities for high-order schemes. On the contrary, using $\tilde c$ these schemes are stable, since the level sets remain equally spaced around the $0$-level set and the representation function is still Lipschitz continuous.  
Looking at Tab. \ref{table:rhombus_tildec_1_Nx}, we can observe that the monotone scheme converges in a lower number of iterations $N_i$ and with lower errors with respect to the correspondent results obtained by using the classical velocity and reported in Tab. \ref{table:rhombus_c_1_Nx}. Note  that we used the same parameters, tolerance and initial datum (the paraboloid, Datum 1) for both velocities to make a fair comparison. 
The AF-LW scheme converges and give always better results in terms of both errors and number of iterations with respect to the monotone scheme with $c$ or $\widetilde c$. 
Also the WENO scheme converges in a lower number of iterations with respect to the monotone one, even if with greater errors. This is due to the fact that the WENO scheme needs to recompute the velocity  $\widetilde c$, which makes the errors accumulate (since we use Runge-Kutta of third order, we need to recompute $\widetilde c$ three times). 
Hence, all the three schemes benefit of the new definition of velocity and the AF-LW scheme with $\widetilde c$ gives the best performances. 
Regarding the speed of the schemes, looking at Tab. \ref{tab:rhombus-CPU} which contains the CPU times in seconds related to the three schemes needed to obtain the results reported in Tab. \ref{table:rhombus_tildec_1_Nx}, the AF-LW scheme needs a longer CPU time with respect to the monotone one, as expected. However, the AF-LW scheme is faster than the  WENO scheme and its  best performances are reached in a short time, at most about 40 seconds for the last refinement of the grid, compared to the $1028.51$ seconds necessary to the WENO scheme.   

\begin{table}[h!]
	\caption{Test 1a. Errors and number of iterations using $L^\infty$ norm, Datum $1$ and the parameters $\mu=2$, $K_{reg}=0$, $tol=0.0005$, varying the number of nodes. Best results are in bold.}
	\centering
	\begin{tabular}{c| c c c | c c c  | c c c  }
		$c$ & \multicolumn{3}{|c|}{\bf Monotone} & \multicolumn{3}{|c}{\bf AF-LW}& \multicolumn{3}{|c}{\bf WENO}\\
		\hline
		$\# Nodes$ &  $N_i$  & $P$-$Err_{rel}$ & $P$-$Err_1$ & $N_i$  & $P$-$Err_{rel}$ & $P$-$Err_1$ & $N_i$  & $P$-$Err_{rel}$ &\hspace{-0.2cm}$P$-$Err_1$\\
		\hline
		$102$  & $\bf 84$	& $\bf 0.1025$ & $\bf 0.2321$ &     $213$	& $X$ & $X$  & $217$	& $X$ &\hspace{-0.2cm}$X$ \\
		\hline
		$202$ &  $\bf 152$	& $\bf 0.0526$ & $\bf 0.1172$ &     $394$	& $X$ & $X$ &  $399$	 & $X$ &\hspace{-0.2cm}$X$	\\  
		\hline
		$402$  &  $\bf 288$	& $\bf 0.0265$ & $\bf 0.0593$ &     $745$	& $X$ & $X$ & $669$	& $X$ &\hspace{-0.2cm}$X$ \\
		\hline
	\end{tabular}
	\label{table:rhombus_c_1_Nx}
\end{table}

\begin{table}[h!]
	\caption{Test 1a. Errors and number of iterations using $L^\infty$ norm, Datum $1$ and the parameters $\mu=2$, $K_{reg}=0$, $tol=0.0005$, varying the number of nodes. Best results are in bold.}
	\centering
	\begin{tabular}{c| c c c | c c c  | c c c  }
		$\widetilde c$ & \multicolumn{3}{|c|}{\bf Monotone} & \multicolumn{3}{|c}{\bf AF-LW}& \multicolumn{3}{|c}{\bf WENO}\\
		\hline
		$\# Nodes$ &  $N_i$  & $P$-$Err_{rel}$ & $P$-$Err_1$ & $N_i$  & $P$-$Err_{rel}$ & $P$-$Err_1$ & $N_i$  & $P$-$Err_{rel}$ &\hspace{-0.2cm}$P$-$Err_1$   \\
		\hline
		$102$  & $50$	& $0.0748$ & $0.1694$ &     $48$	& $\bf 0.0693$ & $\bf 0.1568$  & $\bf 47$	& $0.0886$ &\hspace{-0.2cm}$0.2008$ \\
		\hline
		$202$ &   $100$	& $0.0427$ & $0.0950$ &     $\bf 96$	& $\bf 0.0363$ & $\bf 0.0808$ &  $\bf 96$	 & $0.0469$ &\hspace{-0.2cm}$0.1045$	\\  
		\hline
		$402$  & $199$	& $0.0208$ & $0.0466$ &     $\bf 196$	& $\bf 0.0203$ & $\bf 0.0454$ & $\bf 196$	& $0.0240$ &\hspace{-0.2cm}$0.0537$ \\
		\hline
	\end{tabular}
	\label{table:rhombus_tildec_1_Nx}
\end{table}

\begin{table}[h!]
	\caption{Test 1a. CPU times in seconds related to the Tab. \ref{table:rhombus_tildec_1_Nx}.}
	\centering
	\vspace{0.1cm}
	\begin{tabular}{|c | c | c | c |}
		\hline
		$\# Nodes$ &{\bf Monotone }  & {\bf AF-LW } &{\bf WENO}\\
		\hline
		102 & $0.19$ &  $0.64$  &  $4.75$	\\  
		\hline
		202 & $1.35$ &  $5.61$ &  $68.38$	\\  
		\hline
		402 & $9.77$ &  $40.20$ &  $1028.51$ 	\\  
		\hline
	\end{tabular}
	\label{tab:rhombus-CPU}
\end{table}
A more correct analysis can be made by decreasing the tolerance $tol$ together with the number of nodes $\# Nodes$. In Tabs. \ref{table:rhombus_tildec_1_Datum-monotone}, \ref{table:rhombus_tildec_1_Datum-AF}, and \ref{table:rhombus_tildec_1_Datum-WENO} we analyze the behavior of the schemes with respect to the initial data. Looking at Tab. \ref{table:rhombus_tildec_1_Datum-monotone} we can note that the monotone scheme is not influenced by the change of the initial datum in terms of number of iterations. In fact, comparing the two columns related to $N_i$, only one iteration in the last row is different. 
With respect to the errors, small changes are visible, with a small improvement using Datum 2 for the last two rows. 
The AF-LW scheme (Tab. \ref{table:rhombus_tildec_1_Datum-AF}) has  better performances with lower (or equal) errors and lower (or equal) $N_i$ when using Datum 2. 
For the WENO scheme (see Tab. \ref{table:rhombus_tildec_1_Datum-WENO}) no changes are visible in terms of $N_i$ for all the refinements and both the initial data. 
In terms of errors, only a small change in the last row is visible, with a slight preference for Datum 1. 
For all the schemes, the errors decrease  when we refine the grid as expected.  
Comparing the three tables, Tabs. \ref{table:rhombus_tildec_1_Datum-monotone},  \ref{table:rhombus_tildec_1_Datum-AF}, and \ref{table:rhombus_tildec_1_Datum-WENO}, we can note that the AF-LW scheme get always better results in terms of number of iterations and errors with both initial data with respect to the monotone scheme, expect for the last refinement ($\# Nodes = 402$) with Datum 2, in which the monotone scheme seems to be a little bit better, due to round off errors. 
The WENO scheme get always the greatest errors using both intial data, for the reasons explained before, in a number of iterations which differ from those of the AF scheme for at most one iteration in some cases. 
This simple synthetic example illustrates very well the limits of a high-order scheme like the WENO one, when some singularities occur, and how the proposed AF scheme is able to deal with them. 
\begin{table}[h!]  
	\caption{Test 1a. Errors and number of iterations varying $\# Nodes$ ($L^\infty$ norm). Best results are in bold.}\label{table:rhombus_tildec_1_Datum-monotone}
	\centering
	\begin{tabular}{c|c c c | c c c | c c c  }
		$\widetilde c$ & \multicolumn{3}{|c|} {\bf Monotone} & \multicolumn{3}{|c|}{\bf Datum $1$} & \multicolumn{3}{|c}{\bf Datum $2$}\\
		\hline
		\hline
		$\# Nodes$  & $tol$ &$\mu$  & $K_{reg}$ & $N_i$ & $P$-$Err_{rel}$ & $P$-$Err_1$ & $N_i$ & $P$-$Err_{rel}$ & $P$-$Err_1$   \\
		\hline
		$102$ & $0.001$	& $2$ & $0$ &  $\bf 50$	& $\bf 0.0748$ & $\bf 0.1694$ &     $\bf 50$	& $0.0776$ & $0.1757$     \\
		
		$202$ & $0.0005$	& $2$ & $0$ &  $\bf 100$	& $0.0427$ & $0.0950$ &     $\bf 100$	& $\bf 0.0420$ & $\bf 0.0935$	\\
		
		$402$ & $0.00025$	& $2$ & $0$ &  $\bf 199$	& $0.0208$ & $0.0466$ &      $200$	& $\bf 0.0199$ & $\bf 0.0446$	\\  
		\hline
	\end{tabular}
\end{table}
\begin{table}[h!]  
	\caption{Test 1a. Errors and number of iterations varying $\# Nodes$ ($L^\infty$ norm). Best results are in bold.}\label{table:rhombus_tildec_1_Datum-AF}
	\centering
	\begin{tabular}{c|c c c | c c c | c c c  }
		$\widetilde c$ & \multicolumn{3}{|c|} {\bf  AF-LW} & \multicolumn{3}{|c|}{\bf Datum $1$} & \multicolumn{3}{|c}{\bf Datum $2$}\\
		\hline
		\hline
		$\# Nodes$  & $tol$ &$\mu$  & $K_{reg}$ & $N_i$ & $P$-$Err_{rel}$ & $P$-$Err_1$ & $N_i$ & $P$-$Err_{rel}$ & $P$-$Err_1$   \\
		\hline
		$102$ & $0.001$	& $2$ & $0$ & $48$	& $0.0693$ & $0.1568$   &     $\bf 47$	& $\bf 0.0658$ & $\bf 0.1490$      \\
		
		$202$ & $0.0005$	& $2$ & $0$ &  $\bf 96$	& $\bf 0.0363$ & $\bf 0.0808$  &     $\bf 96$	& $\bf 0.0363$ & $\bf 0.0808$	\\
		
		$402$ & $0.00025$	& $2$ & $0$ &  $196$	& $0.0203$ & $0.0454$ &     $\bf 195$	& $\bf 0.0201$ & $\bf 0.0450$	\\
		\hline
	\end{tabular}
\end{table}
\begin{table}[h!]  
	\caption{Test 1a. Errors and number of iterations varying $\# Nodes$ ($L^\infty$ norm). Best results are in bold.\label{table:rhombus_tildec_1_Datum-WENO}}
	\centering
	\begin{tabular}{c|c c c | c c c | c c c  }
		$\widetilde c$ & \multicolumn{3}{|c|} {\bf  WENO} & \multicolumn{3}{|c|}{\bf Datum $1$} & \multicolumn{3}{|c}{\bf Datum $2$}\\
		\hline
		\hline
		$\# Nodes$  & $tol$ &$\mu$  & $K_{reg}$ & $N_i$ & $P$-$Err_{rel}$ & $P$-$Err_1$ & $N_i$ & $P$-$Err_{rel}$ & $P$-$Err_1$   \\
		\hline
		$102$ & $0.001$	& $2$ & $0$ & $\bf 47$	& $\bf 0.0886$ & $\bf 0.2008$   &     $\bf 47$	& $\bf 0.0886$ & $\bf 0.2008$      \\
		
		$202$ & $0.0005$	& $2$ & $0$ &  $\bf 96$	& $\bf 0.0469$ & $\bf 0.1045$  &     $\bf 96$	& $\bf 0.0469$ & $\bf 0.1045$	\\
		
		$402$ & $0.00025$	& $2$ & $0$ &  $\bf 196$ & $\bf 0.0240$ & $\bf 0.0537$ &     $\bf 196$	& $0.0242$ & $0.0541$	\\
		\hline
	\end{tabular}
\end{table}

\subsection{Real tests}\label{sec:RT}
\begin{figure}[h!]
	\centering
	\includegraphics[width=0.32\textwidth]{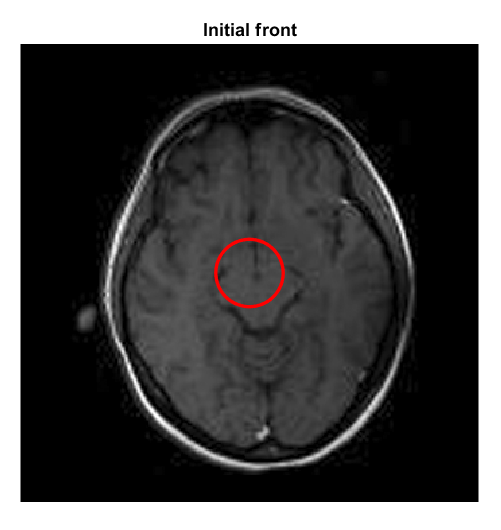}
	\includegraphics[width=0.32\textwidth]{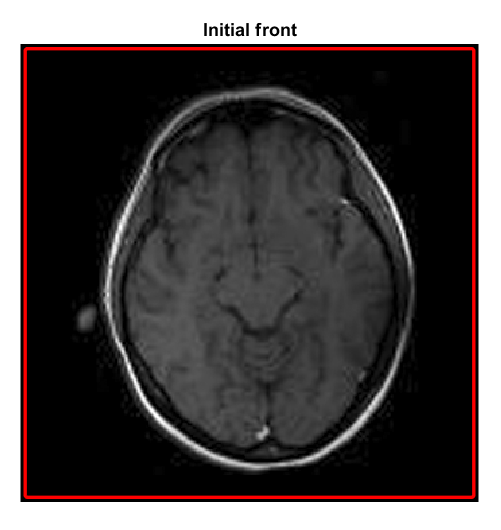}
	\includegraphics[width=0.32\textwidth]{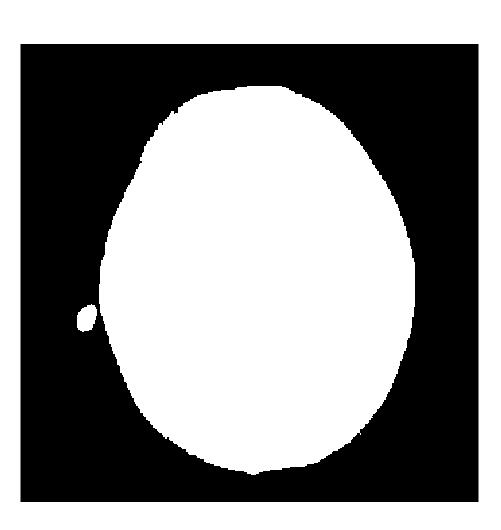}
	\vspace{-0.2 cm}
	\caption{\small{Test 2. From left to right: Initial front for the expansion case (Case 2a) composed by a circle of radius $r=0.25$; initial front for the shrinking case (Case 2b); mask used for the pixel errors in Case 2b.}}
	\label{brain:initial-condition}
\end{figure}

In this section we consider real images also coming from  biomedical applications. 
Thanks to the error formulas \eqref{eq:err-rel-formula} and \eqref{eq:err1-formula}, we can give a sort of quantitative evaluation of the performances of the schemes in terms of ``pixels error" (i.e. the number of pixels composing the area of the object to be segmented). 
\paragraph{Test 2. Brain ($340\times 340$ pixels)}
The first real test focuses on a biomedical image of a human brain. For this test, we approximate the relevant ``external" boundary of the brain via a front expansion (Case a) or a shrinking (Case b), so that we start from inside or outside as visible in Fig. \ref{brain:initial-condition}. 

\begin{figure}[t!]
	\centering
	\includegraphics[width=0.32\textwidth]{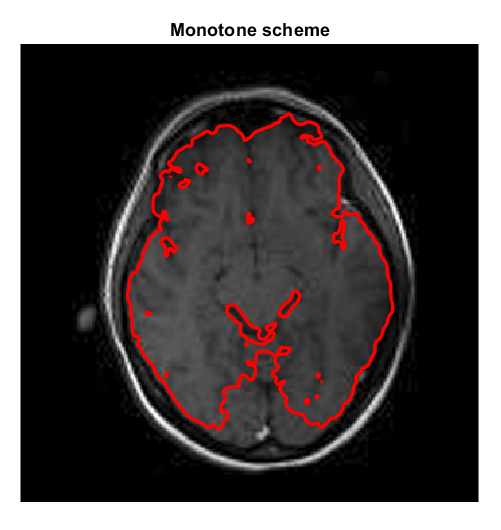}
	\quad
	\includegraphics[width=0.32\textwidth]{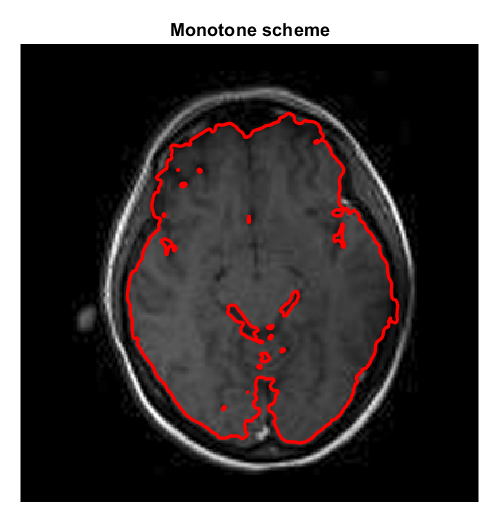}\\
	\includegraphics[width=0.32\textwidth]{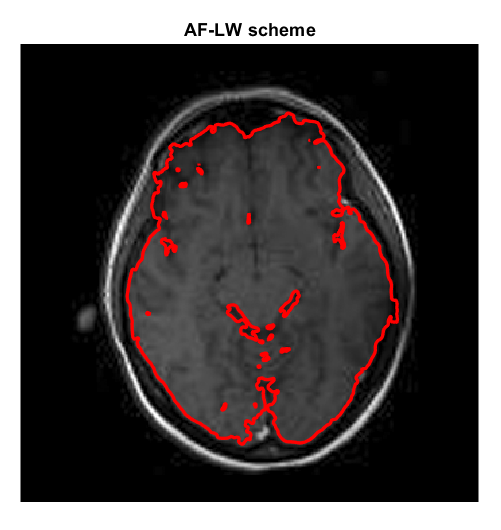}
	\quad
	\includegraphics[width=0.32\textwidth]{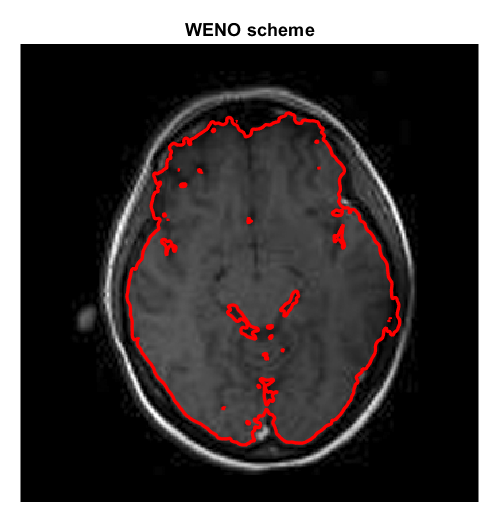}
	\vspace{-0.2 cm}
	\caption{\small{Test 2a with Datum 1. On the first row: plots of the final front using the monotone scheme with $tol=0.00001$, $N_i=376$ (left), and with $tol=0.000005$, $N_i=468$ (right).   
			Second row: plot of the final front using the AF-LW scheme, $N_i=407$ (left), and the WENO scheme, $N_i=451$ (right), both with  $tol=0.00001$. 
			The four tests have been obtained using the $L^1$ norm in the stopping criterion, with $\mu=4$, $K_{reg}=5$, and velocity $\widetilde c$.}}
	\label{brain:exp-datum1}
\end{figure}

\begin{figure}[h!]
	\centering
	\includegraphics[width=0.32\textwidth]{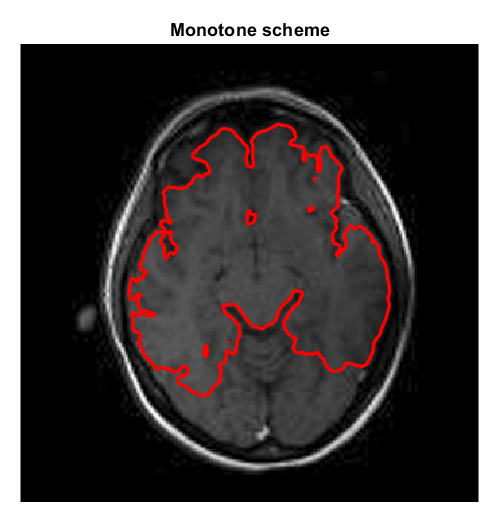}
	\quad
	\includegraphics[width=0.32\textwidth]{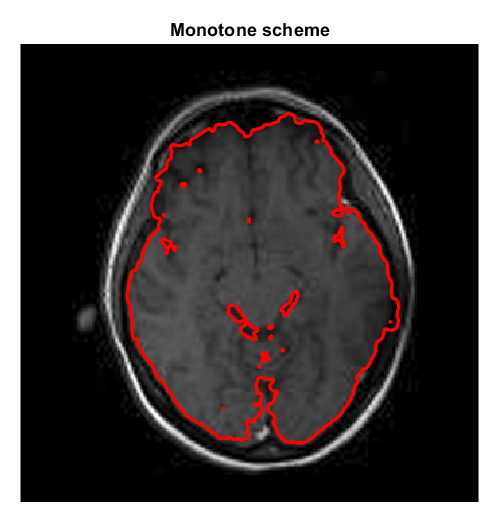}\\
	\includegraphics[width=0.32\textwidth]{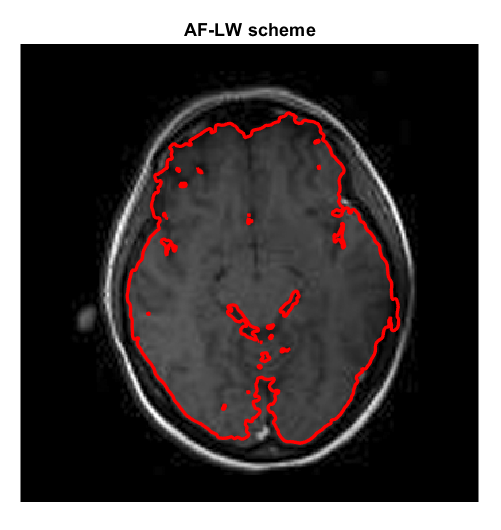}
	\quad
	\includegraphics[width=0.32\textwidth]{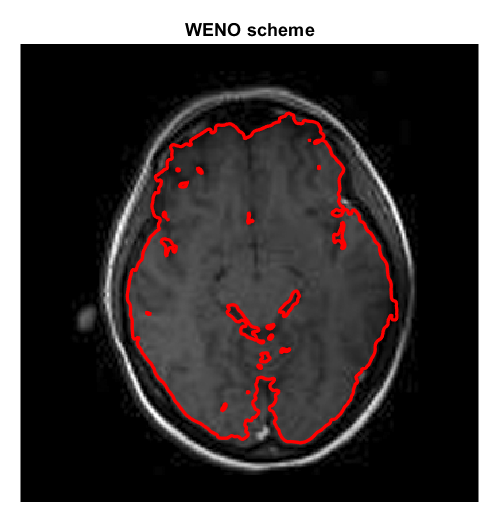}
	\vspace{-0.2 cm}
	\caption{\small{Test 2a with Datum 2. First row: Plots of the final front using the monotone scheme  with $tol=0.00005$, $N_i=264$ (left), and with $tol=0.00001$, $N_i=495$ (right). Second row: Plots of the final front using the AF-LW scheme, $N_i=431$ (left), and the WENO scheme, $N_i=390$ (right), both with  $tol=0.00005$. The four tests have been obtained using the $L^1$ norm in the stopping criterion, with $\mu=4$, $K_{reg}=5$, and velocity $\widetilde c$.}}
	\label{brain:exp-datum2}
\end{figure}

Looking at Figs. \ref{brain:exp-datum1} and \ref{brain:exp-datum2}, it is worth to note that the two high-order schemes recognize better the boundary of the object starting from two different initial datum using the same tolerance $tol=0.00001$. 
The differences and the best resolutions are clearly visible looking at the central part close to the bottom of the brain figures. 
In fact, the monotone scheme stops too early (see the top-left pictures in both Figs. \ref{brain:exp-datum1} and \ref{brain:exp-datum2}). In order to get better results, a smaller tolerance parameter is necessary for the monotone scheme, see the top-right pictures in both the considered figures. 
In that case the monotone scheme can increase its resolution even if with a higher number of iterations and in any case not more accurate than the high-order schemes.  
Comparing the AF and the WENO schemes, both obtain more accurate results with a lower number of iterations with respect to the monotone scheme. More in details, with initial Datum 1, AF scheme converges in a lower number of iterations, with Datum 2 the WENO scheme uses less iterations to get convergence. 
In any case, the CPU times are really different and the AF scheme is always much more faster than the WENO scheme (aroud 15-20 times faster), as visible looking at Tab. \ref{tab:brain-CPU}, in which the CPU times in seconds related to the simulations of the three schemes visible in Figs. \ref{brain:exp-datum1}-\ref{brain:exp-datum2} are reported.   
Clearly, the AF-LW scheme needs  more CPU time with respect to the monotone scheme, this is mainly due to the computation of the smoothness indicators, but  is still fast and competitive since it needs only one minute and half in the worst case (to get a better result). 
In Case a, hence, the qualitative evaluation is enough to understand which scheme provides the better results in a reasonable CPU time. 

\begin{table}[h!]
	\caption{{Test 2a. CPU times in seconds related to the brain tests (Figs. \ref{brain:exp-datum1}-\ref{brain:exp-datum2}).}}
	\centering
	\vspace{0.1cm}
	\begin{tabular}{|c | c | c | c | c |}
		\hline
		Figure &{\bf Mon.(left)} &{\bf Mon.(right) } & {\bf AF-LW }&{\bf WENO}\\
		\hline
		\ref{brain:exp-datum1} & $12.95$	& $16.54$  & $71.18$ & $1424.18$  	\\  
		\hline
		\ref{brain:exp-datum2} & $9.63$ 	& $20.31$  & $91.41$ & $1318.59$	\\  
		\hline
	\end{tabular}
	\label{tab:brain-CPU}
\end{table}

\begin{figure}[h!]
	\centering
	\includegraphics[width=0.32\textwidth]{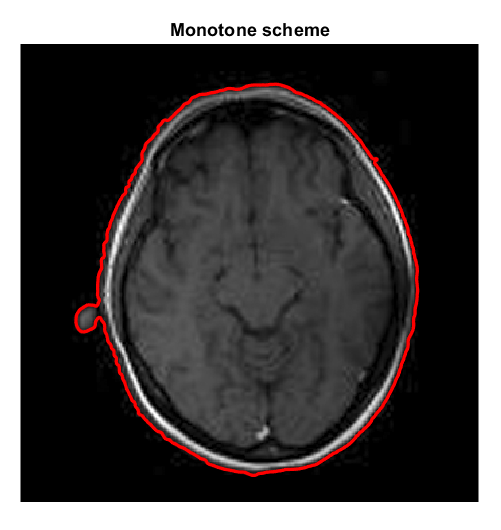}
	\includegraphics[width=0.32\textwidth]{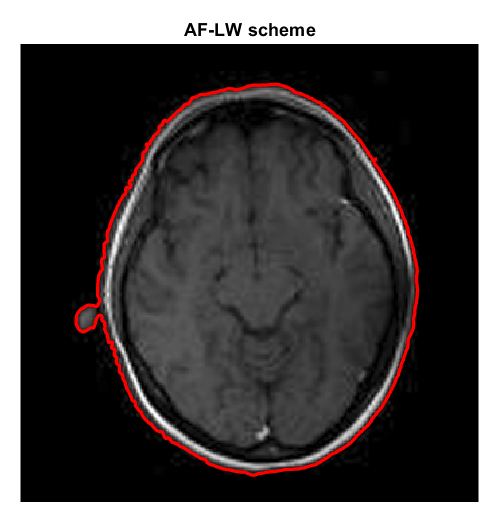}
	\includegraphics[width=0.32\textwidth]{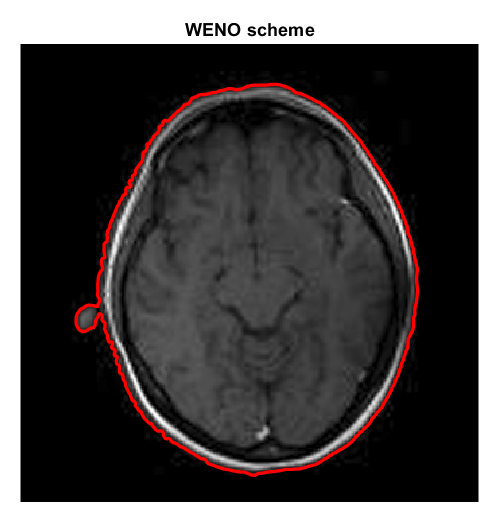}
	\vspace{-0.2 cm}
	\caption{\small{Test 2b. Plots of the final front using the monotone scheme (left), the AF-LW scheme (middle) and the WENO scheme (right) with $\mu=5$, $K_{reg}=3$, and velocity $\widetilde c$. Image size: $340\times340$. }}
	\label{brain:shrink}
\end{figure}

For  the shrinking case, a quantitative error evaluation is needed in addition to the qualitative one. Looking at Fig. \ref{brain:shrink}, we can see that all the three schemes recognized the desired more external boundary.  Analyzing the errors reported in Tab. \ref{tab:brain-shrink}, we observe that the AF-LW scheme produces  lower values in both errors $P$-$Err_{rel}$ and $P$-$Err_1$ with respect to the other schemes for both the resolutions considered (input image size $170\times170$ and $340\times340$ pixels). The errors in the second row have been computed using the mask visible in Fig. \ref{brain:initial-condition} on the right and are related to Fig. \ref{brain:shrink}. In that test of the brain, case b, the WENO scheme obtains errors closer to those of the AF scheme, lower with respect to the monotone scheme for both the considered image sizes. 
Comparing the CPU times with the two different resolutions, we note that 
with half size and less than half time we can obtain better accuracy with respect to the $P$-$Err_1$ error using the AF-LW scheme instead of the monotone one. 
Comparing only the two high-order schemes, again the CPU times of the WENO scheme are about $11$ or $20$ times greater than those necessary to the AF scheme. 
\begin{table}[h!]
	\caption{Test 2b. Errors and number of iterations using $L^1$ norm, $\widetilde c$ and the parameters $tol=0.00005$, $\mu=5$, $K_{reg}=3$, varying the image size. Best results are in bold.}
	\centering
	\begin{tabular}{c| c c c |  c c c  |  c c c }
		&  \multicolumn{3}{|c|}{ \bf Monotone } & \multicolumn{3}{|c}{ \bf AF-LW }& \multicolumn{3}{|c}{\bf WENO}\\
		\hline
		Image size &  $N_i$  & $P$-$Err_{rel}$ &\hspace{-0.2cm}$P$-$Err_1$ & $N_i$  & $P$-$Err_{rel}$ &\hspace{-0.2cm} $P$-$Err_1$ & $N_i$ & $P$-$Err_{rel}$ &\hspace{-0.2cm}$P$-$Err_1$   \\
		\hline
		$170\times170$ &  $\bf 83$ & $0.0465$ &\hspace{-0.2cm}$0.2436$ &  $87$  & $\bf 0.0439$ &\hspace{-0.2cm} $\bf 0.2300$ &  $87$  & $ 0.0448$ &\hspace{-0.2cm}$0.2348$	\\  
		\hline
		$340\times340$  & $\bf 228$ & $0.0118$ &\hspace{-0.2cm}$0.2476$ &  $262$ & $\bf 0.0078$ &\hspace{-0.2cm} $\bf 0.1628$ &  $264$ & $ 0.0079$ &\hspace{-0.2cm}$0.1648$ \\
		\hline
	\end{tabular}
	\label{tab:brain-shrink}
\end{table}

\begin{table}[h!]
	\caption{Test 2b. CPU times in seconds related to the brain tests (Fig. \ref{brain:shrink}).}
	\centering
	\vspace{0.1cm}
	\begin{tabular}{|c | c | c | c |}
		\hline
		Image size &{\bf Monotone } & {\bf AF-LW }&{\bf WENO}\\
		\hline
		$170\times170$ 	& $0.80$  & $3.75$ & $44.10$  	\\  
		\hline
		$340\times340$  	& $8.39$  & $46.07$ & $933.90$	\\  
		\hline
	\end{tabular}
	\label{tab:brain-CPU_case-b}
\end{table}

\begin{figure}[h!]
	\centering
	\includegraphics[width=0.2\textwidth]{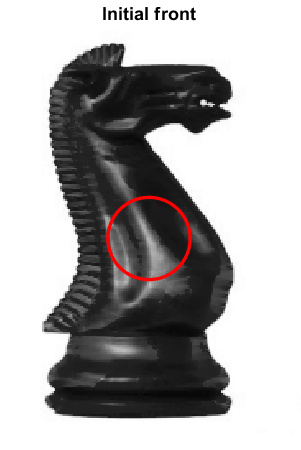}
	\qquad\qquad
	\includegraphics[width=0.2\textwidth]{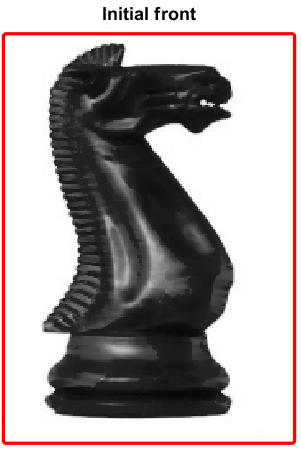}
	\qquad\qquad
	\includegraphics[width=0.2\textwidth]{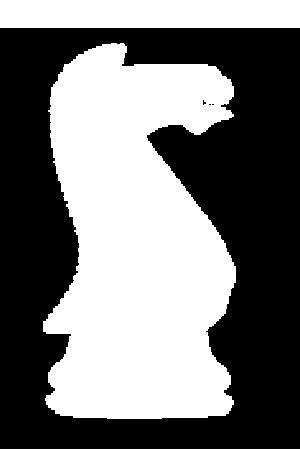}
	\vspace{-0.1 cm}
	\caption{\small{Test 3. From left to right: Initial front for the expansion case (Case 3a) composed by a circle of radius $r=0.5$; initial front for the shrinking case (Case 3b); mask used for the pixel errors in both cases.}}
	\label{horse:initial-condition}
\end{figure}

\begin{figure}[t!]
	\centering
	\includegraphics[width=0.2\textwidth]{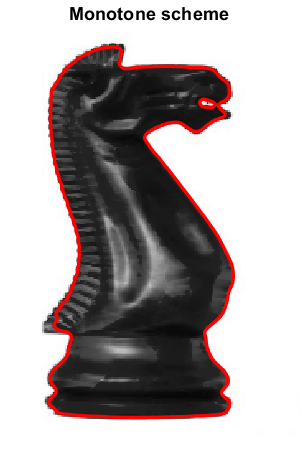}
	\qquad\qquad
	\includegraphics[width=0.2\textwidth]{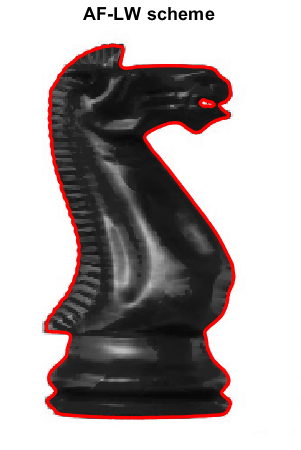}
	\qquad\qquad
	\includegraphics[width=0.2\textwidth]{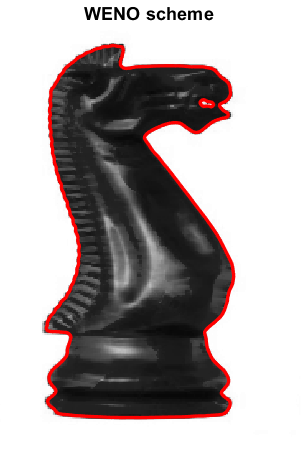}
	\vspace{-0.3 cm}
	\caption{\small{Test 3a with Datum 1. Plots of the final front using the monotone scheme, $N_i=508$ (left), the AF-LW scheme, $N_i=544$ (middle), and the WENO scheme, $N_i=523$ (right), with $L^1$ norm and $tol=0.000025$, $\mu=2$,  $K_{reg}=5$, and velocity $\widetilde c$. 
		}
		\label{horse:expansion_datum1}}
\end{figure}

\begin{figure}[t!]
	\centering
	\includegraphics[width=0.2\textwidth]{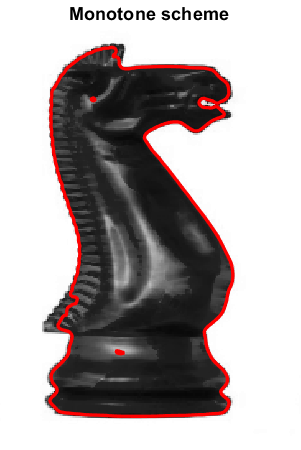}
	\qquad\qquad
	\includegraphics[width=0.2\textwidth]{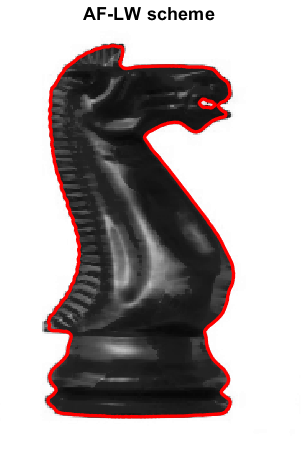}
	\qquad\qquad
	\includegraphics[width=0.2\textwidth]{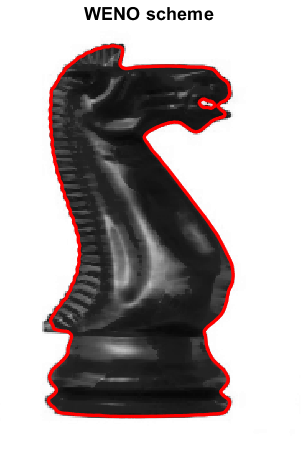}
	\vspace{-0.3 cm}
	\caption{\small{Test 3a with Datum 2. Plots of the final front using the monotone scheme  $N_i=453$ (left), the AF-LW scheme, $N_i=542$ (middle), and the WENO scheme, $N_i=472$ (right),  with $L^1$ norm and $tol=0.00004$, $\mu=2$,  $K_{reg}=5$, and velocity $\widetilde c$. 
		}
		\label{horse:expansion_datum2}}
\end{figure}

\paragraph{Test 3. Horse Chess ($184\times 256$  pixels)}
We choose the image of the horse piece in the game of chess visible in Fig. \ref{horse:initial-condition} and we approximate its boundary from inside (Case a), varying the initial datum, and outside (Case b), always using the modified model with velocity $\widetilde{c}$.  
Starting our analysis of the results from the expansion case, looking at Fig. \ref{horse:expansion_datum1} we can note some differences between the schemes around the mouth of the horse, at the top at the beginning of the horse's mane, and at the bottom left, i.e. at the end of the horse's mane, where the monotone scheme seems to stop too early. 
Looking at Fig. \ref{horse:expansion_datum2}, we note in addition other differences and worse performances of the monotone scheme, due probably to the specularities inside the image. 
Regarding the errors reported in Tab. \ref{tab:horse_a1-err}, the WENO scheme get lower errors with Datum 1 and $tol=0.00005$, even if with a greater number of iterations with respect to the AF scheme, whereas with $tol=0.000025$ we obtain the inverse situation, i.e. the AF scheme get lower errors but with more iterations with respect to the WENO scheme. 
Looking at Tab. \ref{tab:horse_a2-err}, we can note that using the initial Datum 2 the AF scheme provides always the best performances in terms of both errors $P$-$Err_{rel}$ and $P$-$Err_1$.
Regarding the CPU times, Tab. \ref{tab:horse_a-CPU} shows that the WENO scheme is really slow compared to the other schemes, whereas the AF scheme always converges in less than one minute. 
\begin{table}[h!]
	\caption{Test 3a. Errors and number of iterations using $L^1$ norm, Datum $1$, velocity $\widetilde c$, and the parameters $\mu=2$, $K_{reg}=5$, varying the tolerance. Best results are in bold.}
	\centering
	\begin{tabular}{c| c c c | c c c  | c c c  }
		$\widetilde c$ & \multicolumn{3}{|c|}{\bf Monotone} & \multicolumn{3}{|c}{\bf AF-LW}& \multicolumn{3}{|c}{\bf WENO}\\
		\hline
		tol &  $N_i$  & $P$-$Err_{rel}$ & $P$-$Err_1$ & $N_i$  & $P$-$Err_{rel}$ &\hspace{-0.2cm} $P$-$Err_1$ & $N_i$  & $P$-$Err_{rel}$ &\hspace{-0.2cm}$P$-$Err_1$   \\
		\hline
		$0.00005$ &  $433$	 & $0.0794$ & $0.6424$ &  $\bf 396$	 & $0.0672$ &\hspace{-0.2cm} $0.6424$ &  $432$	 & $\bf 0.0591$ &\hspace{-0.2cm}$\bf 0.4788$	\\  
		\hline
		$0.000025$  & $\bf 508$	 & $0.0599$ & $0.4848$ &  $544$	& $\bf 0.0438$ &\hspace{-0.2cm} $\bf 0.3544$  & $523$	& $0.0462$ &\hspace{-0.2cm}$0.3744$ \\
		\hline
	\end{tabular}
	\label{tab:horse_a1-err}
\end{table}
\begin{table}[h!]
	\caption{Test 3a. Errors and number of iterations using $L^1$ norm, Datum $2$, velocity $\widetilde c$, and the parameters $\mu=2$, $K_{reg}=5$, varying the tolerance. Best results are in bold.}
	\centering
	\begin{tabular}{c| c c c | c c c  | c c c  }
		$\widetilde c$ & \multicolumn{3}{|c|}{\bf Monotone} & \multicolumn{3}{|c}{\bf AF-LW}& \multicolumn{3}{|c}{\bf WENO}\\
		\hline
		tol &  $N_i$  & $P$-$Err_{rel}$ & $P$-$Err_1$ & $N_i$  & $P$-$Err_{rel}$ & $P$-$Err_1$ & $N_i$  & $P$-$Err_{rel}$ & $P$-$Err_1$   \\
		\hline
		$0.00008$  & $\bf 382$	 & $0.1031$ & $0.8344$ &  $451$	& $\bf 0.0640$ & $\bf 0.5180$  & $411$	& $0.0684$ & $0.5536$ \\
		\hline
		$0.00004$ &  $\bf 453$	 & $0.0737$ & $0.5964$ &  $542$	 & $\bf 0.0501$ & $\bf 0.4052$ &  $472$	 & $0.0565$ & $0.4572$	\\  
		\hline
	\end{tabular}
	\label{tab:horse_a2-err}
\end{table}
\begin{table}[h!]
	\caption{Test 3a. CPU times in seconds related to the chess horse tests (Tabs. \ref{tab:horse_a1-err}-\ref{tab:horse_a2-err}).}
	\centering
	\vspace{0.1cm}
	\begin{tabular}{|c | c | c | c | c |}
		\hline
		Datum & tol &{\bf Monotone }  & {\bf AF-LW } &{\bf WENO}\\
		\hline
		$1$ & $0.00005$ & $7.26$ 	& $26.22$  & $378.07$ \\
		\hline
		$1$ & $0.000025$ & $8.78$ 	& $42.23$  & $462.55$ 	\\  
		\hline
		$2$ & $0.00008$ & $6.51$ 	& $37.15$  & $377.47$ 	\\ 
		\hline
		$2$ & $0.00004$ & $8.02$ 	& $44.16$  & $443.66$ 	\\ 
		\hline
	\end{tabular}
	\label{tab:horse_a-CPU}
\end{table}
\begin{figure}[h!]
	\centering
	\includegraphics[width=0.2\textwidth]{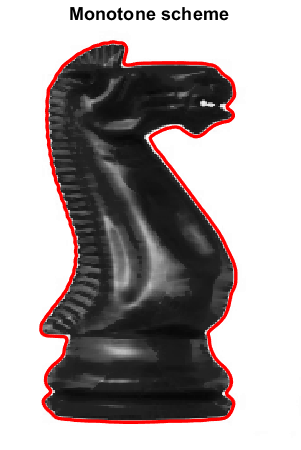}
	\qquad\qquad
	\includegraphics[width=0.2\textwidth]{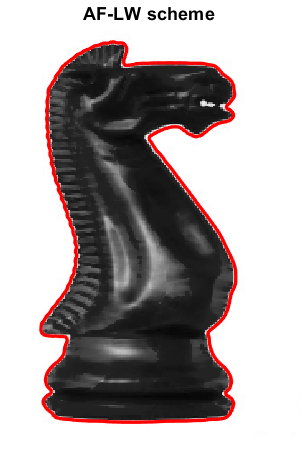}
	\qquad\qquad
	\includegraphics[width=0.2\textwidth]{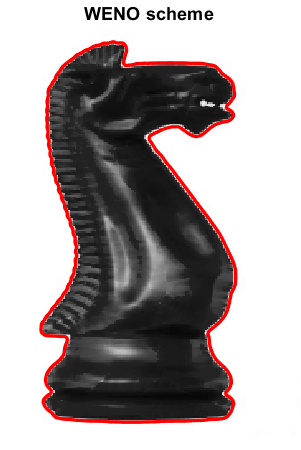}
	\vspace{-0.3 cm}
	\caption{\small{Test 3b. Plots of the final front using the monotone scheme (left), the AF-LW scheme (middle) and the WENO scheme (right) with velocity $\widetilde c$, and $tol=0.00005$, by using $L^1$ norm in the stopping criterion and parameters $\mu=2$, and $K_{reg}=3$. 
		}
		\label{fig:horse_shrink}}
\end{figure}
\begin{table}[h!]
	\caption{Test 3b. Errors and number of iterations using $L^1$ norm, velocity $\widetilde c$ and the parameters $tol=0.00005$, $\mu=2$, $K_{reg}=3$. Best  results are in bold.}
	\centering
	\begin{tabular}{c c c  | c c c  | c c c  }
		\multicolumn{3}{c|}{\bf Monotone} & \multicolumn{3}{|c}{\bf AF-LW}& \multicolumn{3}{|c}{\bf WENO}\\
		\hline
		$N_i$  & $P$-$Err_{rel}$ & $P$-$Err_1$ & $N_i$  & $P$-$Err_{rel}$ & $P$-$Err_1$ & $N_i$  & $P$-$Err_{rel}$ & $P$-$Err_1$   \\
		\hline
		$\bf 187$  & $0.0303$ & $0.2452$ &  $197$   & $0.0273$ & $0.2212$ &  $196$   & $\bf 0.0262$ & $\bf 0.2120$	\\  
		\hline
	\end{tabular}
	\label{tab:horse_b-err}
\end{table}
\begin{table}[h!]
	\caption{Test 3b. CPU times in seconds related to the horse test (Tab. \ref{tab:horse_b-err}).}
	\centering
	\vspace{0.1cm}
	\begin{tabular}{| c | c | c |}
		\hline
		{\bf Monotone }  & {\bf AF-LW } &{\bf WENO}\\
		\hline
		$3.14$ 	& $14.10$  & $195.13$ 	\\      
		\hline
	\end{tabular}
	\label{tab:horse_b-CPU}
\end{table}

\begin{figure}[h!]
	\centering
	\includegraphics[width=0.3\textwidth]{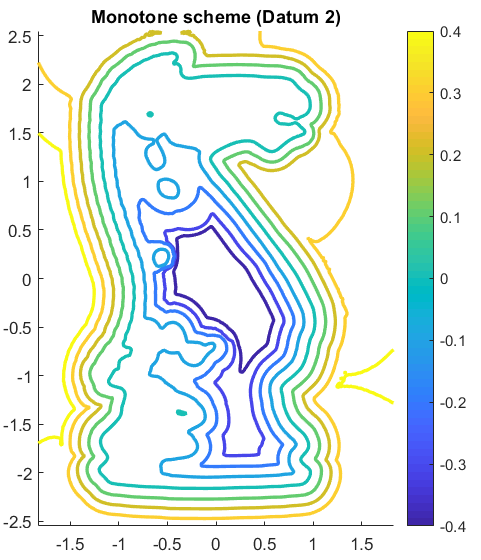}
	\includegraphics[width=0.3\textwidth]{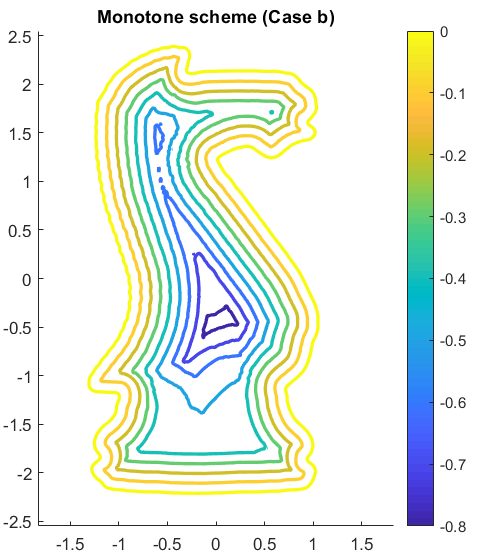}
	\includegraphics[width=0.3\textwidth]{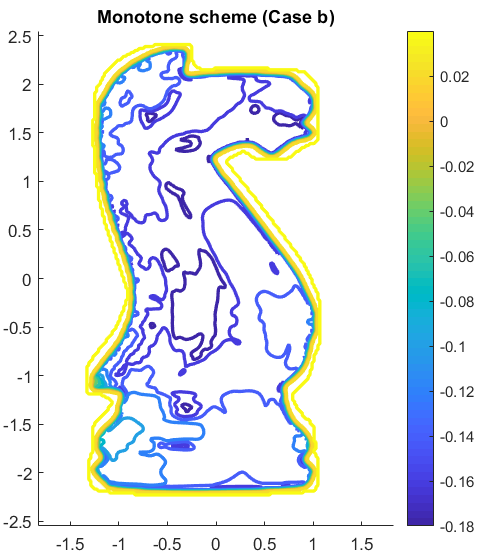}
	\caption{\small{Test 3. Contour plots of the final representations using the monotone scheme with velocity $\widetilde{c}$, for Case a with Datum 2 (left) and for Case b (middle), and with velocity $c$ for Case b (right). }
		\label{fig:horse_repr_monotone}}
\end{figure}

\noindent Analyzing the results in the shrinking case, small differences can be noted looking at Fig. \ref{fig:horse_shrink}, particularly between the two high-order schemes. For the monotone scheme, we note that the final front stops just before the boundary of the object. So, also in this case, the only qualitative analysis is not enough. A quantitative analysis is shown in Tab. \ref{tab:horse_b-err}, in which both the high-order schemes are more accurate than the monotone scheme as expected.  In this case WENO is the most accurate scheme, having lower values in both errors, even if with a CPU time around 14 times greater than that of the AF scheme (see Tab. \ref{tab:horse_b-CPU}). \\
In order to show the effectiveness of our implementation of the modified velocity, in Fig. \ref{fig:horse_repr_monotone} we collected the final representations obtained by the monotone scheme. 
Although in the first case some new fronts arise due to the specularities inside the image, we can still recognize that the gradient of the initial condition is preserved and that the level sets do not collide during the evolution. In particular, in the middle image the distance function to the $0$-level set is still clearly visible. The third image is related to the results obtained by using the classical velocity $c$ in Case b, with the same tolerance used with $\widetilde{c}$ and visible in Fig. \ref{fig:horse_shrink}. We note that with the classical model oscillations appear throughout the whole horse, even if focusing on the $0$-level set the scheme performs well. A better behavior of all the level-sets can be obtained but using a greater tolerance $tol$, e.g. $tol= 0.0005$. 

\begin{figure}[h!]
	\centering
	\includegraphics[width=0.32\textwidth]{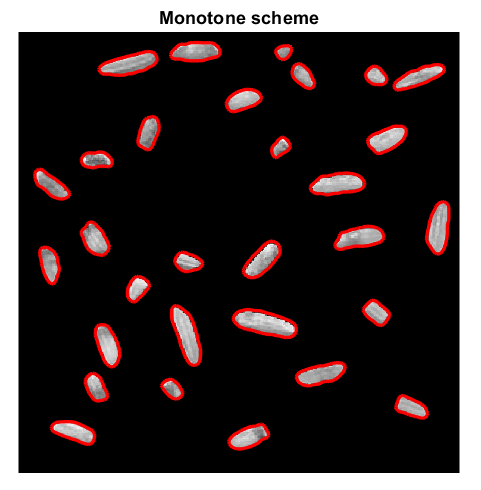}
	\includegraphics[width=0.32\textwidth]{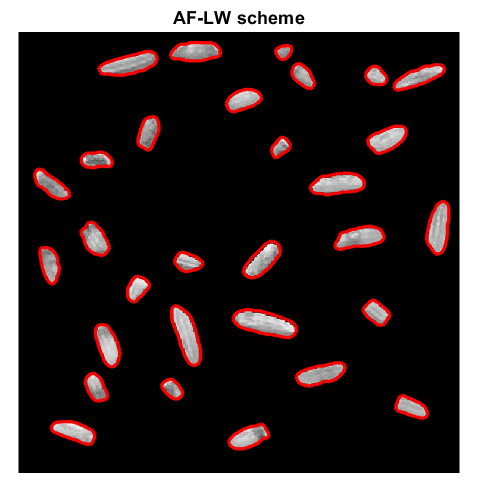}
	\includegraphics[width=0.32\textwidth]{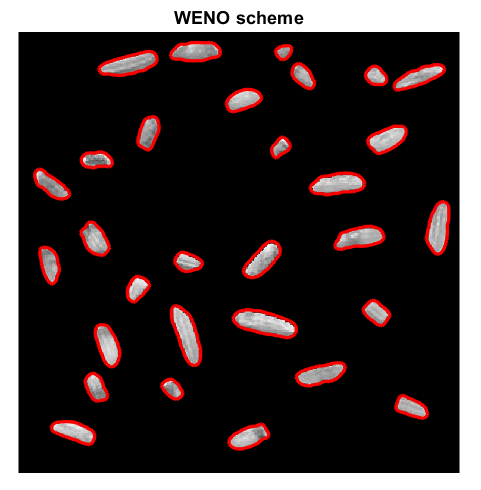}
	\vspace{-0.2 cm}
	\caption{\small{Test 4b. Plots of the final front obtained by the Monotone scheme (left), the AF-LW scheme (middle) and the WENO scheme (right), using velocity $\widetilde{c}$ and the parameters reported in Tab. \ref{tab:grains-err}.}
		\label{fig:fronts-grains}}
\end{figure}

\begin{figure}[h!]
	\centering
	\includegraphics[width=0.43\textwidth]{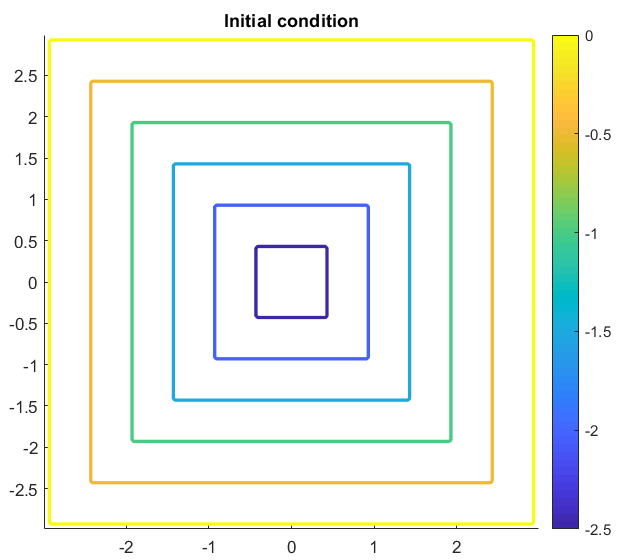}\qquad
	\includegraphics[width=0.44\textwidth]{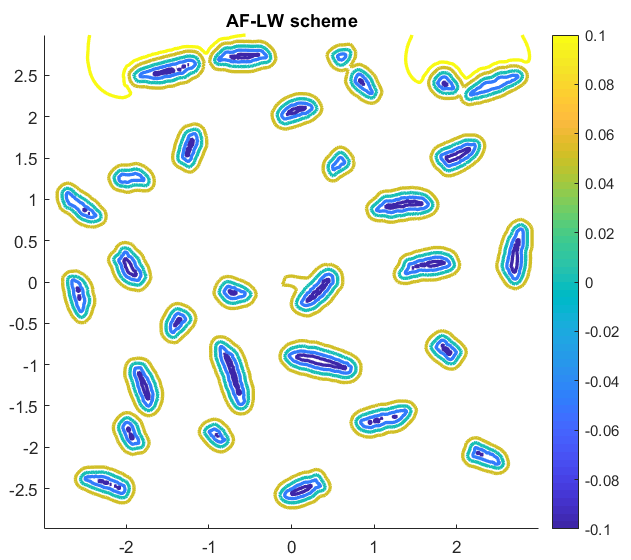}
	\vspace{-0.2 cm}
	\caption{\small{Test 4b. Contour plots of the initial datum (left) and the final representations (right), using the AF-LW scheme with  velocity $\widetilde{c}$  and the parameters reported in Tab. \ref{tab:grains-err}.}
		\label{fig:repr-grains}}
\end{figure}

\paragraph{Test 4. Grains ($300\times 300$ pixels)}
In this test we consider the shrinking front in presence of multiple separate objects to be segmented.  
The final fronts obtained by the three schemes are visible in Fig. \ref{fig:fronts-grains}, in which almost no differences are visible. 
Looking at the errors and number of iterations reported in Tab. \ref{tab:grains-err}, we can note that the two high-order schemes are more accurate than the first order monotone one, with a number of iterations $N_i$ very close between the three schemes, the lowest number is obtained by the AF scheme. In terms of errors, the WENO scheme is the most accurate, even if with a really great CPU time as ever (see Tab. \ref{tab:grains-CPU}).
We report in Fig. \ref{fig:repr-grains} the contour plots of the initial datum and the final representation obtained by the AF scheme. It is rather interesting to see that, although the representation function has a rather complex evolution and splits in various parts, the final solution is still a distance function to the $0$-level set, now composed by different peaks each relative to a single grain. 
\begin{table}[h!]
	\caption{Test 4b. Errors and number of iterations with velocity $\widetilde c$, using $L^1$ norm in the stopping criterion, and the parameters $tol=0.0001$, $\mu=2$, and $K_{reg}=2$. Best results are in bold.}
	\centering
	\begin{tabular}{c c c | c c c  | c c c  }
		\multicolumn{3}{c|}{\bf Monotone} & \multicolumn{3}{|c}{\bf AF-LW}& \multicolumn{3}{|c}{\bf WENO}\\
		\hline
		$N_i$  & $P$-$Err_{rel}$ & $P$-$Err_1$ & $N_i$  & $P$-$Err_{rel}$ & $P$-$Err_1$ & $N_i$  & $P$-$Err_{rel}$ & $P$-$Err_1$   \\
		\hline
		$312$	& $0.0097$ & $0.0316$ &     $\bf 299$	& $0.0074$ & $0.0240$ & $306$	& $\bf 0.0064$ & $\bf 0.0208$	\\  
		\hline
	\end{tabular}
	\label{tab:grains-err}
\end{table}
\begin{table}[h!]
	\caption{Test 4b. CPU times in seconds related to the grains test (Tab. \ref{tab:grains-err}).}
	\centering
	\vspace{0.1cm}
	\begin{tabular}{|c| c| c|}
		\hline
		{\bf Monotone }  & {\bf AF-LW } &{\bf WENO}\\
		\hline
		$8.71$ 	& $41.26$  & $676.91$ 	\\   
		\hline
	\end{tabular}
	\label{tab:grains-CPU}
\end{table}

\begin{figure}[h!]
	\centering
	\includegraphics[width=0.65\textwidth]{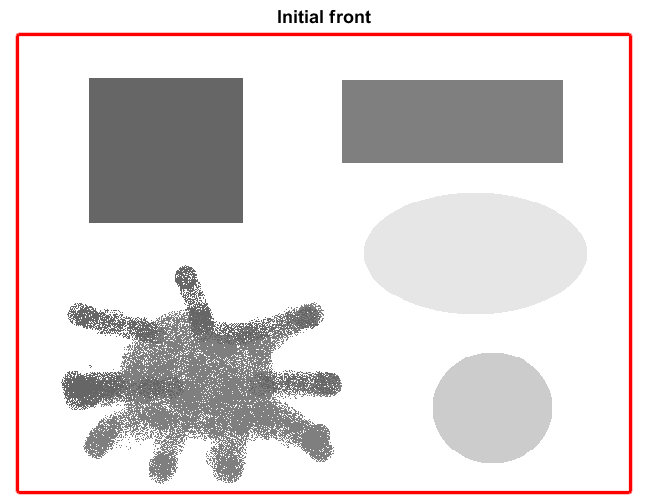}
	\caption{\small{Test 5b. Initial front for the shrinking case (Case 5b).}}
	\label{shapes:initial_data}
\end{figure}

\begin{figure}[h!]
	\centering
	\includegraphics[width=0.44\textwidth]{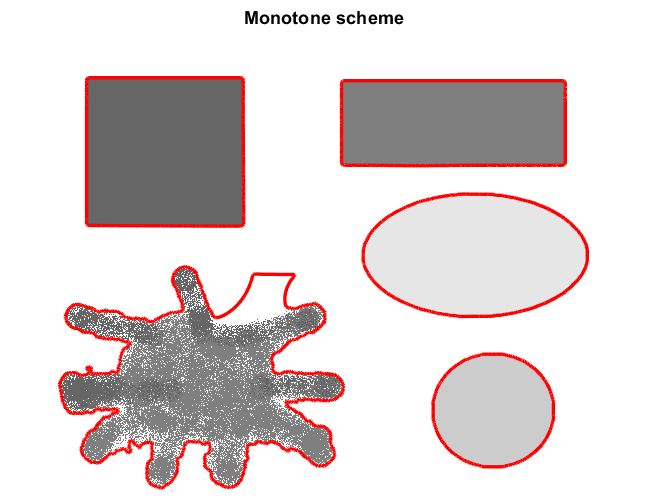}
	\includegraphics[width=0.44\textwidth]{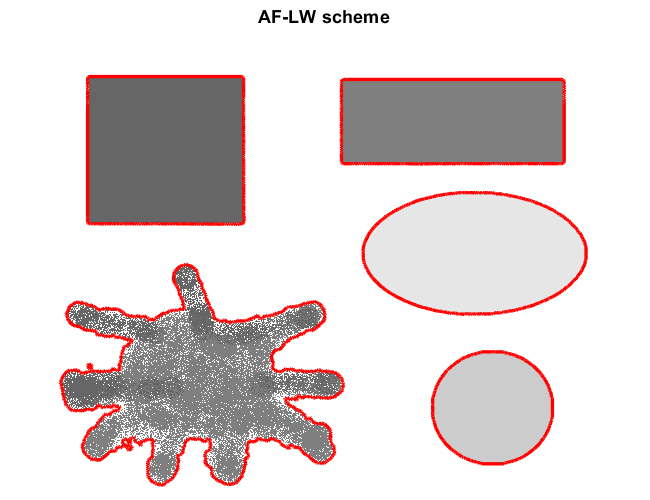}\\
	\quad
	\includegraphics[width=0.44\textwidth]{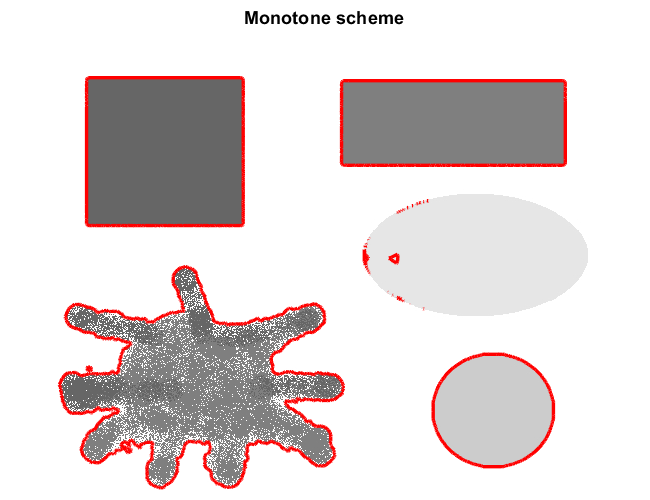}
	\includegraphics[width=0.44\textwidth]{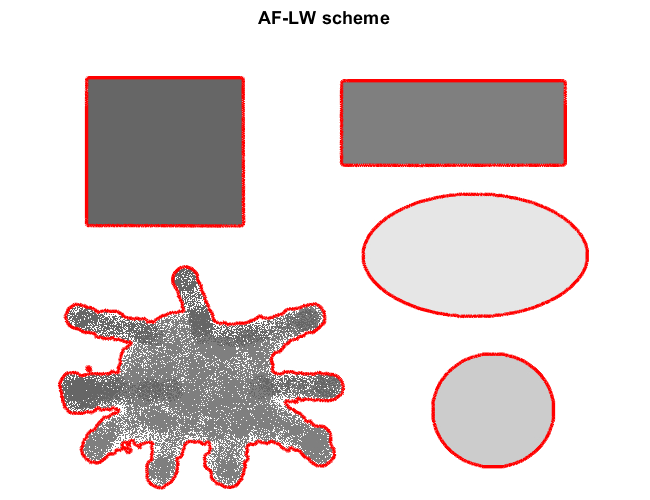}\\
	\quad
	\includegraphics[width=0.44\textwidth]{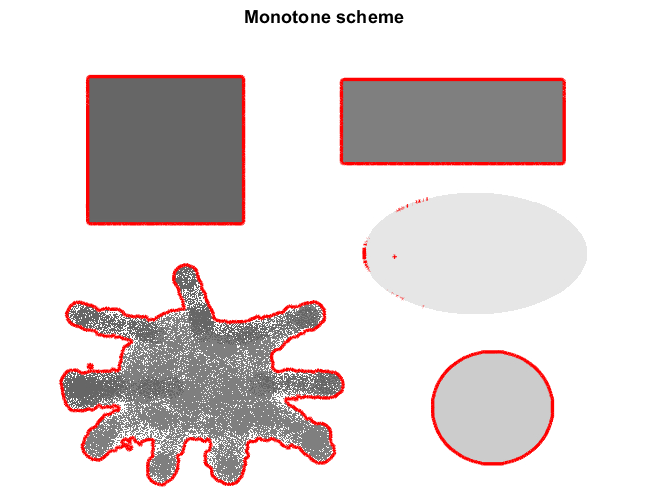}
	\includegraphics[width=0.44\textwidth]{Forme_AFLW_p4N1_tol0_00005_arresto2.png}
	\vspace{-0.2 cm}
	\caption{\small{Test 5b. Plots of the final front using the monotone scheme (left), and the AF-LW scheme varying the tolerance ($tol=0.0005,0.0001,0.00005$) with $\mu=4$, $K_{reg}=1$, and velocity $\widetilde c$. 
	}}
	\label{shapes:shrink}
\end{figure}

\paragraph{Test 5. Geometric shapes ($640\times480$ pixels)}
In some of the previous real tests, we have seen that with a careful tuning of the parameters in the stopping rule, the monotone scheme can get comparable results with respect to those obtained by the high-order schemes, at least from a visible point of view. This is not always possible, especially in more critical situations, e.g. when the difference between the background and the objects we want to segment is less marked, as shown in this test. 
In Fig.  \ref{shapes:initial_data} we reported the initial front in red for the shrinking case. The more difficult object to detect will be the ellipse, due to the lighter gray levels closer to the white background.  
Looking at the final fronts in Fig. \ref{shapes:shrink}, obtained by varying the tolerance in the stopping criterion, we can note that the AF scheme recognizes very well the boundaries of all the objects, differently from the monotone one. In fact, if we stop too early the schemes, with $tol=0.0005$, we note that the monotone scheme still has to conclude the recognition of the ``spray" shape below on the left. But if we adopt a smaller tolerance, e.g. $tol=0.0001$ or $tol=0.00005$, in order to give ``more time" to the monotone scheme to achieve all the boundaries, the scheme improves the detection of the spray shape, but loses the boundary of the ellipse. 
Higher-order approximation of the evolution can clearly give more stability, especially when the contrast is not satisfactory as in the present situation. \\
This numerical test clearly shows the difficulties that the monotone scheme can encounter. It needs a really complicated manual tuning of the parameters, and nonetheless not always gives good and reliable results (as in this case), whereas the AF scheme overcomes this difficulty thanks to its properties. 

\begin{figure}[h!]
	\centering
	\includegraphics[width=0.325\textwidth]{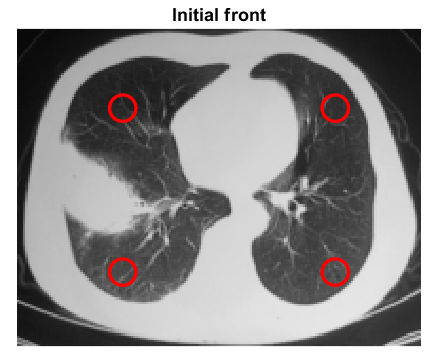}
	\includegraphics[width=0.325\textwidth]{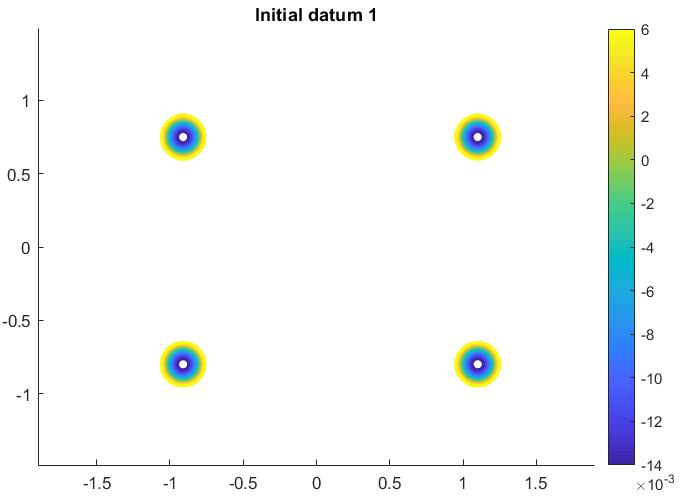}
	\includegraphics[width=0.325\textwidth]{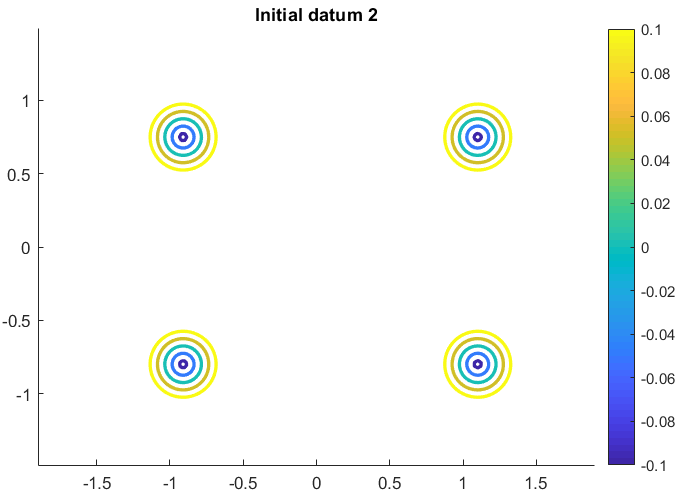}
	\caption{\small{Test 6a. From left to right: Initial front composed of four separate circles of radius $r=0.125$, contour plots of the initiam data (Datum 1, and Datum 2).}
		\label{fig:pneumonia_initial}}
\end{figure}

\paragraph{Test 6. Pneumonia ($191\times150$ pixels)}
Finally, we conclude our numerical tests considering a Pneumonia image in the expansion case starting from an initial datum here composed by four equal paraboloids (each defined as Datum 1) or cones (defined as Datum 2), placed as visible in Fig. \ref{fig:pneumonia_initial}, with $0$-level sets composed by circles of radius $r=0.125$. 

\begin{figure}[h!]
	\centering
	\includegraphics[width=0.32\textwidth]{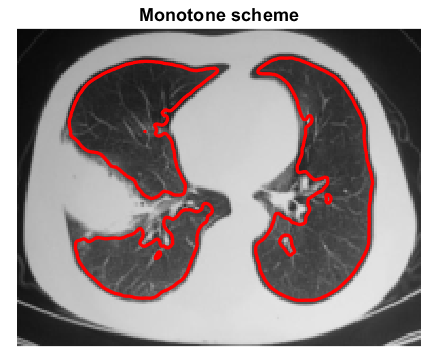}
	\includegraphics[width=0.32\textwidth]{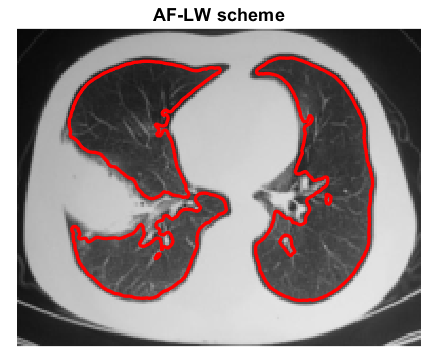}
	\includegraphics[width=0.32\textwidth]{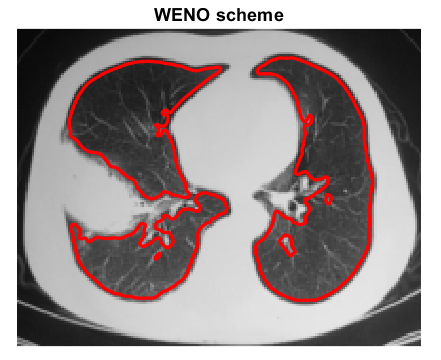}
	\vspace{-0.3cm}
	\caption{\small{Test 6a with Datum 1. From left to right: Plots of the final front using the monotone scheme ($N_i=185$), the AF-LW scheme ($N_i=194$), and the WENO scheme ($N_i=192$), with $L^1$ norm and $tol=0.00001$, $\mu=4$, $K_{reg}=5$, and velocity $\widetilde c$.}
		\label{fig:pneumonia_datum1}}
\end{figure}

\begin{figure}[h!]
	\centering
	\includegraphics[width=0.32\textwidth]{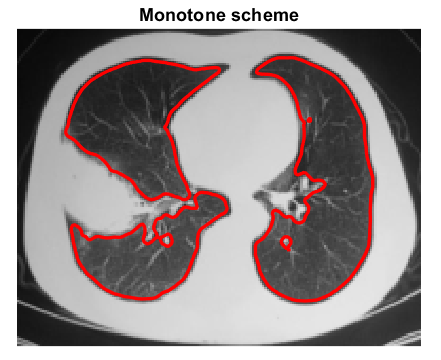}
	\includegraphics[width=0.32\textwidth]{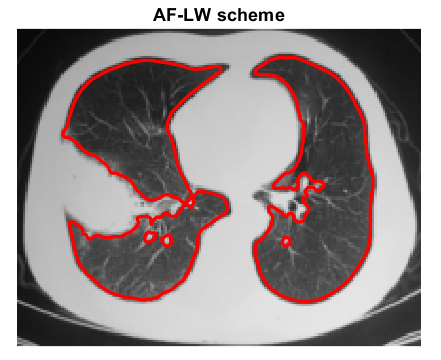}
	\includegraphics[width=0.32\textwidth]{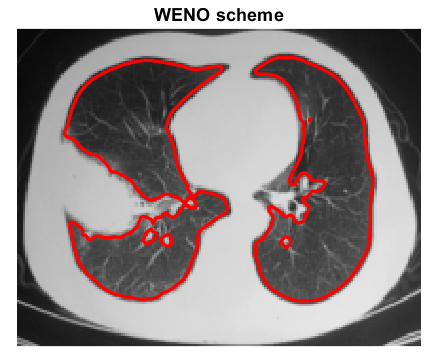}
	\vspace{-0.3cm}
	\caption{\small{Test 6a with Datum 2. From left to right: Plots of the final front using the monotone scheme ($N_i=413$), $tol=0.00001$, the AF-LW scheme ($N_i=589$) with $tol=0.00001$, and the WENO scheme ($N_i=458$) with $tol=0.000012$, all using $L^1$ norm, $\mu=4$, $K_{reg}=5$, and velocity $\widetilde c$.}
		\label{fig:pneumonia_datum2}}
\end{figure}

\noindent In this test the behavior of the three schemes is confirmed in terms of CPU time, looking at Tab. \ref{tab:pneumonia-CPU}, and in terms of accuracy, as visible in  Figs. \ref{fig:pneumonia_datum1} and \ref{fig:pneumonia_datum2}, where the final fronts obtained by the high-order schemes recognize better the object starting from the two different initial data, Datum 1 and Datum 2, respectively. In fact, in the first figure, Fig. \ref{fig:pneumonia_datum1}, it is evident looking at the lung on the left;  
in the second case, Fig. \ref{fig:pneumonia_datum2}, it is clear especially looking at the right lung, mostly for visualizing the better performances of the AF scheme. 
Moreover, comparing in vertical the two figures, focusing on each scheme with different initial datum, we can note that all the three schemes seems to prefer the distance function (Datum 2). 

\begin{figure}[h!]
	\centering
	\includegraphics[width=0.325\textwidth]{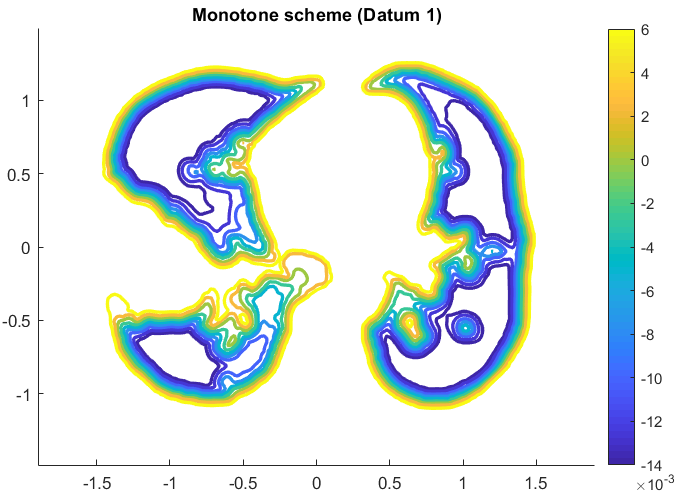}
	\includegraphics[width=0.325\textwidth]{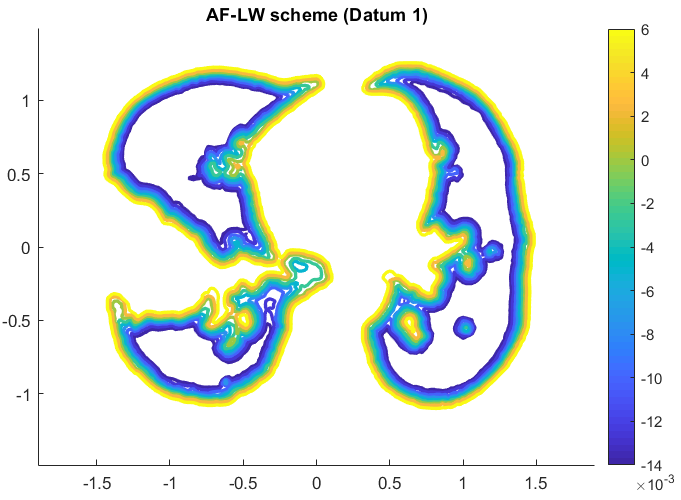}
	\includegraphics[width=0.325\textwidth]{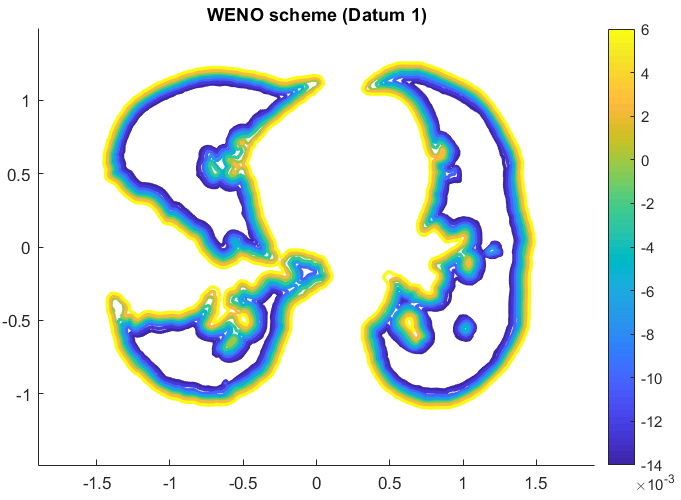}\\
	\includegraphics[width=0.325\textwidth]{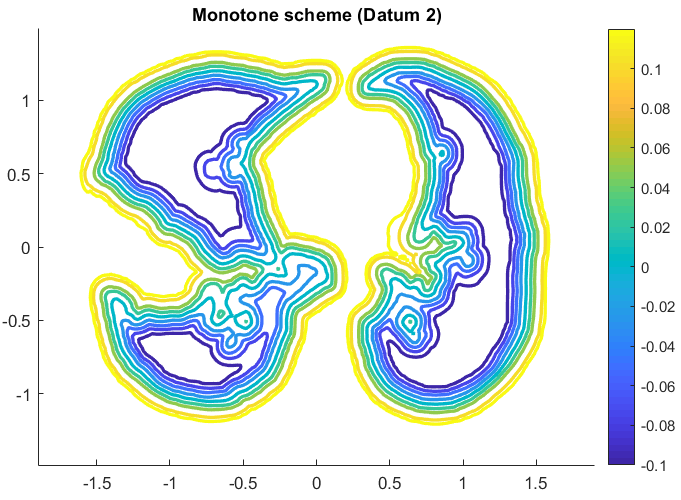}
	\includegraphics[width=0.325\textwidth]{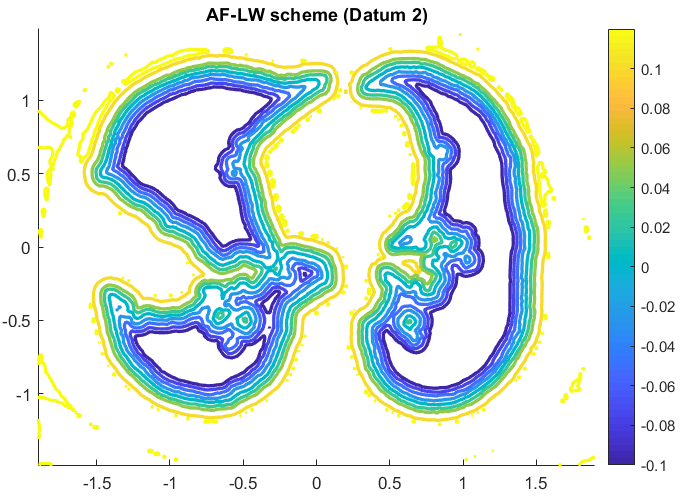}
	\includegraphics[width=0.325\textwidth]{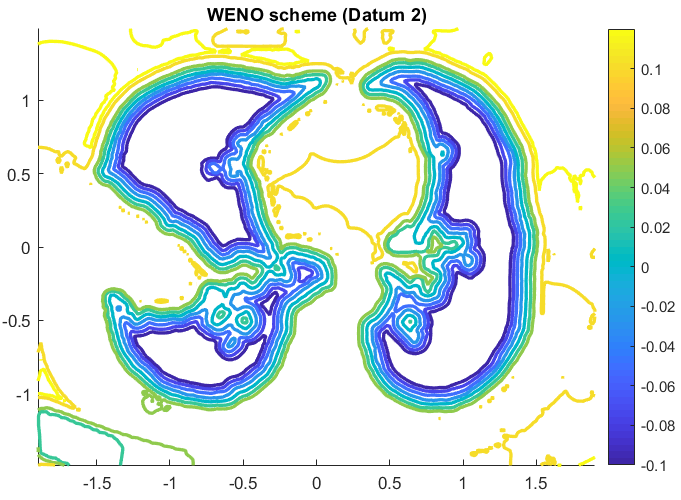}
	\caption{\small{Test 6a. Contour plots of the final representations using the monotone scheme (left), the AF-LW scheme (middle), and the WENO scheme (right). Top: Datum 1, bottom: Datum 2.}
		\label{fig:pneumonia_repr_schemes}}
\end{figure}

\noindent The good behavior of the AF scheme and of the modified model in the case of a merging fronts is confirmed by the contour plots of the final representations, shown in Fig. \ref{fig:pneumonia_repr_schemes}. Regarding the WENO scheme (last column),  some oscillations arise at some point in the evolution and increase as time flows (see the bottom-left part of the final representation using Datum 2). This is why a slightly greater tolerance has been used in that case to obtain the reported result.  
\begin{table}[h!]
	\caption{Test 6a. CPU times in seconds related to the Pneumonia tests (Figs. \ref{fig:pneumonia_datum1}-\ref{fig:pneumonia_datum2}).}
	\centering
	\vspace{0.1cm}
	\begin{tabular}{|c|c| c| c|}
		\hline
		Figure &{\bf Monotone}& {\bf AF-LW}& {\bf WENO}\\
		\hline
		\ref{fig:pneumonia_datum1} & $1.80$ &  $8.05$  &  $83.83$	\\  
		\hline
		\ref{fig:pneumonia_datum2} & $5.01$   & $28.91$  & $202.373$ 	\\ 
		\hline
	\end{tabular}
	\label{tab:pneumonia-CPU}
\end{table}

\section{Conclusions and perspectives}\label{sec:conclusions}
In this work we have proposed a new velocity function for  the level-set method in order to improve  image segmentation  and we have extended to 2D an adaptive filtered scheme originally developed in 1D \cite{FPT18a}.   We have shown that the use of the new velocity function $\widetilde c$  allows to get more accurate results stabilizing high-order schemes as the AF or the WENO schemes, for which the classical velocity introduces some instabilities destroying the convergence. Moreover, the new velocity function can also be applied to  the simple monotone scheme getting better results even in the first order approximation and does not require the re-initialization of the front in order to keep an accurate tracking of the 0-level set.
From the numerical point of view, the AF scheme is based on two building blocks: a monotone scheme and a high-order scheme.
The filter function allows to couple the two schemes in a rather simple way and easily switches from one scheme to the other according to new smoothness indicators.
In terms of CPU time, the AF scheme is less expensive than the WENO scheme and offers a good option to improve the accuracy of the monotone scheme.
It is interesting to note that these two changes are rather effective for the segmentation of synthetic and real images according to many simulations, some of them are presented in Sect. \ref{sec:tests}. 
A qualitative analysis of the resulting segmentations is illustrated by the pictures of the last section whereas a more accurate comparison of the schemes is based on a quantitative analysis of the pixel errors.   Moreover, the sensitivity of the AF scheme seems to be rather low with respect to the presence of noise and only  few steps of a linear filter are required to obtain the necessary regularization of the gray levels of the input image $I$.
Although in this paper the image segmentation is obtained just using the classical first order equation, a possible extension to second order problems can be considered, this extension is motivated by the inclusion of curvature terms in the evolutive equation as done in the literature. The analysis of the adaptive filtered scheme to second order non linear equations  goes beyond the scopes of this paper and will be object of a further investigation. 

\section*{Acknowledgments}
We would like to thank the National Group INdAM-GNCS for the financial support given to this research. 


\begin{thebibliography}{10}
	
	\bibitem{AACD02}
	S.~Amat, F.~Ar\`andiga, A.~Cohen, and R.~Donat.
	\newblock {Tensor product multiresolution analysis with error control for
		compact image representation}.
	\newblock {\em Signal Processing}, 82(4):587--608, 2002.
	
	\bibitem{ABM10}
	F.~Ar{\`a}ndiga, A.~M. Belda, and P.~Mulet.
	\newblock {Point-Value WENO Multiresolution Applications to Stable Image
		Compression}.
	\newblock {\em Journal of Scientific Computing}, 43(2):158--182, May 2010.
	
	\bibitem{AB03}
	F.~Ar\`andiga and A.M. Belda.
	\newblock {Weighted ENO interpolation and applications}.
	\newblock {\em Commun. Nonlinear Sci. Numer. Simul.}, 9(2):187--195, 2003.
	
	\bibitem{Barles94}
	G.~Barles.
	\newblock {\em {Solutions de viscosit\`e des equations de Hamilton-Jacobi}}.
	\newblock Springer Verlag, 1994.
	
	\bibitem{BFS16}
	O.~Bokanowski, M.~Falcone, and S.~Sahu.
	\newblock {An efficient filtered scheme for some first order
		Hamilton-Jacobi-Bellman equations}.
	\newblock {\em SIAM J. Sci. Comput.}, 38(1):A171--A195, 2016.
	
	\bibitem{BPR18}
	O.~Bokanowski, A.~Picarelli, and C.~Reisinger.
	\newblock {High-order filtered schemes for time-dependent second order HJB
		equations}.
	\newblock {\em ESAIM Math. Model. Numer. Anal.}, 52(1):69--97, 2018.
	
	\bibitem{BW04}
	T.~Brox and J.~Weickert.
	\newblock {Level Set Based Image Segmentation with Multiple Regions}.
	\newblock In Carl~Edward Rasmussen, Heinrich~H. B{\"u}lthoff, Bernhard
	Sch{\"o}lkopf, and Martin~A. Giese, editors, {\em Pattern Recognition}, pages
	415--423, Berlin, Heidelberg, 2004. Springer Berlin Heidelberg.
	
	\bibitem{CCZ13}
	X.~Cai, R.~Chan, and T.~Zeng.
	\newblock {A Two-Stage Image Segmentation Method Using a Convex Variant of the
		Mumford--Shah Model and Thresholding}.
	\newblock {\em SIAM Journal on Imaging Sciences}, 6(1):368--390, 2013.
	
	\bibitem{CCF06}
	E.~Carlini, E.~Cristiani, and N.~Forcadel.
	\newblock {A non-monotone fast marching scheme for a Hamilton-Jacobi equation
		modeling dislocation dynamics}.
	\newblock In A.~B. de~Castro, D.~G\`omez, P.~Quintela, and P.~Salgado, editors,
	{\em Numerical mathematics and advanced applications, Proceedings of ENUMATH
		2005}, pages 723--731. Springer, 2006.
	
	\bibitem{CFF10}
	E.~Carlini, M.~Falcone, and R.~Ferretti.
	\newblock {Convergence of a large time-step scheme for mean curvature motion}.
	\newblock {\em Interfaces and Free Boundaries}, 12(4):409--441, 2010.
	
	\bibitem{CFF19}
	E.~Carlini, M.~Falcone, and R.~Ferretti.
	\newblock {Numerical techniques for level set models: an image segmentation
		perspective}.
	\newblock In A.~El-Baz and J.~Suri, editors, {\em Level set methods in medical
		imaging segmentation}. Taylor \& Francis, 2019.
	
	\bibitem{CFFM08}
	E.~Carlini, M.~Falcone, N.~Forcadel, and R.~Monneau.
	\newblock {Convergence of a generalized fast-marching method for an eikonal
		equation with a velocity-changing sign}.
	\newblock {\em SIAM J. Numerical Analysis}, 46(6):2920--2952, 2008.
	
	\bibitem{CCCD93}
	V.~Caselles, F.~Catt\'e, T.~Coll, and F.~Dibos.
	\newblock {A geometric model for active contours in image processing}.
	\newblock {\em Num. Math.}, 66:1--31, 1993.
	
	\bibitem{CKS97}
	V.~Caselles, R.~Kimmel, and G.~Sapiro.
	\newblock {Geodesic Active Contours}.
	\newblock {\em International Journal of Computer Vision}, 22(1):61--79, Feb
	1997.
	
	\bibitem{CC91}
	L.D. Cohen.
	\newblock {On active contour models and balloons}.
	\newblock {\em {CVGIP:} Image understanding}, 53(2):211--218, 1991.
	
	\bibitem{CC92}
	L.D. Cohen and I.~Cohen.
	\newblock {Deformable models for 3-D medical images using finite elements and
		balloons}.
	\newblock In {\em Proceedings 1992 IEEE Computer Society Conference on Computer
		Vision and Pattern Recognition}, volume~28, pages 592--598, 1992.
	
	\bibitem{CL83}
	M.G. Crandall and P.-L. Lions.
	\newblock {Viscosity solutions of Hamilton-Jacobi equations}.
	\newblock {\em Transactions of the American Mathematical Society},
	277(1):1--42, 1983.
	
	\bibitem{FPT18a}
	M.~Falcone, G.~Paolucci, and S.~Tozza.
	\newblock {Convergence of Adaptive Filtered Schemes for First Order
		Evolutionary Hamilton-Jacobi Equations}.
	\newblock 2018, submitted.
	\newblock arXiv:1812.02140.
	
	\bibitem{FPT19}
	M.~Falcone, G.~Paolucci, and S.~Tozza.
	\newblock {Adaptive Filtered Schemes for first order Hamilton-Jacobi
		equations}.
	\newblock In F.A. Radu, K.~Kumar, I.~Berre, J.M. Nordbotten, and I.S. Pop,
	editors, {\em Numerical Mathematics and Advanced Applications ENUMATH 2017},
	pages 389--398, Cham, 2019. Springer International Publishing.
	
	\bibitem{FO13}
	B.D. Froese and A.M. Oberman.
	\newblock {Convergent filtered schemes for the Monge-Amp\'ere partial
		differential equation}.
	\newblock {\em SIAM Journal on Numerical Analysis}, 51(1):423--444, 2013.
	
	\bibitem{HEOC86}
	A.~Harten, S.~Osher, B.~Engquist, and S.~R. Chakravarthy.
	\newblock {Some results on uniformly high-order accurate essentially
		nonoscillatory schemes}.
	\newblock {\em Appl. Numer. Math.}, 2(3-5):347--377, 1986.
	
	\bibitem{HAP05}
	A.~K. Henrick, T.~D. Aslam, and J.~M. Powers.
	\newblock {Mapped weighted essentially non-oscillatory schemes: Achieving
		optimal order near critical points}.
	\newblock {\em J. Comput. Phys.}, 207:542--567, 2005.
	
	\bibitem{JP00}
	G.~Jiang and D.-P. Peng.
	\newblock {Weighted ENO schemes for Hamilton-Jacobi equations}.
	\newblock {\em SIAM Journal on Scientific Computing}, 21:2126--2143, 2000.
	
	\bibitem{JS96}
	G.-S. Jiang and C.-W. Shu.
	\newblock { Efficient implementation of weighted ENO schemes}.
	\newblock {\em J. Comput. Phys.}, 126(1):202--228, 1996.
	
	\bibitem{KKOTY96}
	S.~Kichenassamy, A.~Kumar, P.~Olver, A.~Tannenbaum, and A.~Yezzi.
	\newblock {Conformal curvature flows: from phase transitions to active vision}.
	\newblock {\em Archive for Rational Mechanics and Analysis}, 134(3):275--301,
	1996.
	
	\bibitem{LS95}
	P.L. Lions and P.~Souganidis.
	\newblock {Convergence of MUSCL and filtered schemes for scalar conservation
		laws and Hamilton--Jacobi equations}.
	\newblock {\em Num. Math.}, 69:441--470, 1995.
	
	\bibitem{LOC94}
	X.D. Liu, S.~Osher, and T.~Chan.
	\newblock {Weighted essentially non-oscillatory schemes}.
	\newblock {\em J. Comput. Phys.}, 115(1):200--212, 1994.
	
	\bibitem{MSV95}
	R.~Malladi, J.~A. Sethian, and B.~C. Vemuri.
	\newblock {Shape modeling with front propagation: a level set approach}.
	\newblock {\em IEEE Transactions on Pattern Analysis and Machine Intelligence},
	17(2):158--175, Feb 1995.
	
	\bibitem{OS15}
	A.~M. Oberman and T.~Salvador.
	\newblock {Filtered schemes for Hamilton-Jacobi equations: a simple
		construction of convergent accurate difference schemes}.
	\newblock {\em Journal of Computational Physics}, 284:367--388, 2015.
	
	\bibitem{OF03}
	S.~Osher and R.~Fedkiw.
	\newblock {\em {Level Set Methods and Dynamic Implicit Surfaces}}.
	\newblock Springer, 2003.
	
	\bibitem{OS88}
	S.~Osher and J.~A. Sethian.
	\newblock {Fronts propagating with curvature-dependent speed: Algorithms based
		on Hamilton-Jacobi formulations}.
	\newblock {\em Journal of Computational Physics}, 79(1):12--49, November 1988.
	
	\bibitem{OS91}
	S.~Osher and C.-W. Shu.
	\newblock {High-order essentially non oscillatory schemes for Hamilton-Jacobi
		equations}.
	\newblock {\em SIAM J. Numer. Anal.}, 28(4):907--922, 1991.
	
	\bibitem{PaolucciPhD}
	G.~Paolucci.
	\newblock {\em {Adaptive Filtered Schemes for first order Hamilton-Jacobi
			equations and applications}}.
	\newblock PhD thesis, Dipartimento di Matematica, Sapienza - Universit\`a di
	Roma, Italy, July 2018.
	
	\bibitem{SR09}
	B.~Scheuermann and B.~Rosenhahn.
	\newblock {Analysis of Numerical Methods for Level Set Based Image
		Segmentation}.
	\newblock In G.~Bebis et~al. (eds), editor, {\em Advances in Visual Computing},
	pages 196--207, Berlin, Heidelberg, 2009. Springer Berlin Heidelberg.
	
	\bibitem{Sethian85}
	J.~A. Sethian.
	\newblock Curvature and the evolution of fronts.
	\newblock {\em Commun. in Mathematical Physics}, 101:487--499, 1985.
	
	\bibitem{Sethian99}
	J.A. Sethian.
	\newblock {\em {Level Set Methods and Fast Marching Methods: Evolving
			Interfaces in Computational Geometry, Fluid Mechanics, Computer Vision, and
			Materials Science}}.
	\newblock Cambridge University Press, 2nd edition, 1999.
	
	\bibitem{ZZLZ16}
	K.~Zhang, L.~Zhang, K.-M. Lam, and D.~Zhang.
	\newblock {A level set approach to image segmentation with intensity
		inhomogeneity}.
	\newblock {\em IEEE Transactions on Cybernetics}, 46(2):546--557, 2016.
	
	\bibitem{ZZSZ10}
	K.~Zhang, L.~Zhang, H.~Song, and W.~Zhou.
	\newblock {Active contours with selective local and global segmentation: a new
		formulation and level set method}.
	\newblock {\em Image and Vision Computing}, 28:668--676, 2010.
	
	\bibitem{ZS16}
	Y.T. Zhang and C.-W. Shu.
	\newblock {ENO and WENO schemes}.
	\newblock In {\em Handbook of numerical methods for hyperbolic problems}, pages
	103--122. Elsevier/North-Holland, 2016.
	
	\bibitem{ZMCL16}
	H.~Zhu, J.~Meng, J.~Cai, and S.~Lu.
	\newblock {Beyond pixels : A comprehensive survey from bottom-up to semantic
		image segmentation and cosegmentation}.
	\newblock {\em Journal of Visual Communication and Image Representation},
	34:12--27, 2016.
	
\end{thebibliography}

\end{document}